\documentclass[12pt]{article}

\usepackage[dvips]{epsfig}
\usepackage{mathrsfs}

\usepackage{amsmath,latexsym,amssymb,amsfonts,amsbsy}
\usepackage{bm} 

\numberwithin{equation}{section}
\newtheorem{lemma}{Lemma}[section]
\newtheorem{theorem}{Theorem}[section]
\newtheorem{assumption}{Assumption}[section]

\newtheorem{example}{\sc Example}[section]

\def\proclaim#1{\par \smallskip\noindent {\bf #1}\bgroup\it\ }
\def\endproclaim{\egroup\par\smallskip}

\def\pr{\textsf{P}} 
\def\ep{\textsf{E}} 
\def\Cov{\textsf{Cov}} 
\def\Var{\textsf{Var}} 
\def\Cal#1{{\mathcal #1}}
\def\bk#1{{\bm #1}}

\def\text#1{\mbox{\rm #1}}
\def\overset#1#2{\stackrel{#1}{#2} }

\begin{document}

\title
{\Large \bf \boldmath Adaptive Allocation Theory in Clinical
Trials\\
--A Review  \footnote{ Project supported by  the National
Natural
Science Foundation of China  (No. 10471126 ).}}%

\author{\large Li-Xin Zhang
    \\ \normalsize\emph{Department of Mathematics, Zhejiang University,}
    \\ \normalsize\emph{Hangzhoou 310027, China.}
    \\ \normalsize\emph{E-mail: stazlx@zju.edu.cn}
}

\date{\vspace{-12mm}}

\maketitle

\begin{abstract}
 Various adaptive randomization procedures (adaptive designs)
have been proposed to clinical  trials. This paper discusses several
broad families of procedures, such as the play-the-winner rule and
Markov chain model, randomized play-the-winner rule and urn models,
drop-the-loser rule, doubly biased coin adaptive design. Asymptotic
theories are presented with several pivotal proofs. The effect of
delayed responses, the power and variability comparison of these
designs are also discussed.
\end{abstract}

\section{Introduction.}
\setcounter{equation}{0}

As reported by the World Health Organization (Global Summary of the
AIDS Epidemic, December 2006), the estimated number of people living
with HIV is 39.5 million, causing 2.9 million deaths in 2006 and
13\%  are children under 15 years. The alarming magnitude of AIDS
epidemic and outbreaks of other fetal contagious diseases such as
SARS reveal how vulnerable our health care system is.  In order to
search for more effective treatments, efficient clinical studies are
urgently needed. In clinical trials, the traditional balanced (or
50\%-50\%) treatment allocation rule
 been challenged due to its possible unethical consequences.
A frequently quoted clinical trial is the study of the drug AZT in
reducing risk of maternal-infract HIV transmission. While half of
the pregnant women (239) are given the AZT drug, the remaining
mothers (238) receive the placebo when 50\%-50\% allocation scheme
is used. Only 20 infants are HIV-positive in AZT group and 60 in the
placebo group. (c.f., Connor et al., New England J. Medicine, 1994).
Balanced allocation resulted in many failures in the placebo group.

  Yao and Wei (1996) redesigned the AZT trial
using an adaptive allocation rule, the randomized play-the-wiener
rule proposed by Wei and Durham (1978), and showed a reduction of
several treatment failures under adaptive allocation.
  Adaptive designs, an important subdivision of experimental designs nowadays,
are   allocation rules in which the probability a treatment assigned
to the coming patient depends upon the results of the previous
patients in the study. The basic goal is to skew allocation
probabilities to favor better treatment performance.

Early important work on adaptive designs was carried out by
 Thompson (1933) and Robbins (1952).  Since then, a
steady stream of research (Zelen (1969), Wei and Durham (1978), Wei
(1979), Eisele and Woodroofe (1995), etc) in this area has generated
various treatment allocation schemes for clinical trials. This paper
provides the recent theories of several broad families of designs.
In Section \ref{section2}, we state some limit results on
martingale, which are the basic tools to derive the asymptotic
properties of adaptive designs. In Section \ref{section3}, we
consider Zelen's play-the-winner rule and its generations by Lin, et
al (2003). In Section \ref{section4}, we derive the asymptotic
properties of the adaptive designs based on urn models, a large
family of randomization procedures. In Section \ref{sectionad4}, the
drop-the-loser rule is introduced. In Section \ref{section5}, an
important family of target-driven designs, the doubly adaptive
biased coin designs, are discussed. In Section \ref{section7}, the
effect of the delay of treatment results is discussed. In Section \ref{section6},
we compare the variabilities of different type of adaptive designs.
The lower bound of the asymptotic variability for a pre-specified
allocation proportion is established, and asymptotic best adaptive
designs are provided. Finally, further discussion and future topics
are mentioned in Section \ref{section8}. For the convenience of
reading, we also give several pivotal proofs. The principle ideas of
deriving asymptotic properties can be found from these proofs. Those
who are not interested in the theoretical results can skip these
proofs and go quickly to the last two sections.

The following notations and definitions are introduced to describe
the randomized treatment allocation schemes. Given a clinical trial
with $K$ treatments. Let $\bm X_1, \bm X_2,...$ be the sequence of
random treatment assignments.  For the $m$-th subject,
 $\bm
X_m=(X_{m,1},\ldots,X_{m,K})$ represents the assignment of treatment
such that if the $m$-th subject is allocated to treatment $k$, then
all elements in $\bm X_m$ are $0$ except for the $k$-th component,
$X_{m,k}$, which is $1$. Let $N_{n,k}$ be the number of subjects
assigned to treatment $k$ in the first $n$ assignments and write
$\bm N_n=(N_{n,1},\ldots, N_{n,K})$. Then $\bm N_n=\sum_{m=1}^n\bm
X_m$. We are interested in the statistical behavior of proportions
$N_{n,k}/n$, $k=1,\ldots, K. $

\section{ Preliminaries, limit  theorems on martingales.}\label{section2}
\setcounter{equation}{0}

 The martingale approach is the basic tool to
investigate the asymptotic properties of adaptive designs.
  In this section, we state some limit theorems on martingales.
  For more results, one can refer to Hall and Heyde (1980) and Stout (1974) or other text books.
  Let $\{M_n, \Cal F_n; n\ge 1\}$ be a real martingale sequence
with $\Delta M_n=M_n-M_{n-1}$ being it difference.

\begin{proclaim}{Theorem A} ({\em LLN}) Let $\eta_n>0$ be a sequence of random variables such $\eta_n$ is
$\Cal F_{n-1}$ measurable and $\eta_n\nearrow $ a.s., and $p$ is a
real number in $(0,2]$.

(a) Then with probability one, $M_n/\eta_n\to 0$ on the event $\{
\eta_n \to \infty, \; \sum_{m=1}^{\infty} \ep\big[|\Delta
M_m|^p|\Cal F_{m-1}\big]/\eta_m^p<\infty\}$.

(b) If $\sum_{m=1}^n\ep[ (\Delta M_m)^2|\Cal F_{m-1}]\le \eta_n$
a.s., then with probability one $M_n/\eta_n\to 0$  on the event
$\{\eta_n\to \infty\}$.
\end{proclaim}
\begin{proclaim}{Theorem B} ({\em LIL}) Suppose $\sup_m\ep[|\Delta M_m|^p]<\infty $  for some $p>2$, then
$ M_n=O(\sqrt{n\log\log n})$ a.s.
\end{proclaim}
\begin{proclaim}{Theorem C} ({\em CLT}) Suppose there is a constant $\sigma\ge 0$ such that
$\sum_{m=1}^n \ep[(\Delta M_m)^2|\Cal F_{m-1}]/n \overset{P}\to \sigma^2$. Further assume that
the conditional Lindberg condition
$$ \frac{1}{n}\sum_{m=1}^n \ep\big[(\Delta M_m)^2I\{|\Delta M_m|\ge \epsilon \sqrt{n}| \Cal F_{m-1}\big]\overset{P}\to 0,
\;\; \forall \epsilon>0$$
is satisfied. Then $M_n/\sqrt{n}\overset{\mathscr D}\to N(0,\sigma^2)$. Further,
$M_{[nt]}/\sqrt{n}\overset{\mathscr D}\to \sigma W(t)$ in $D[0,\infty)$, where $W(t)$ is a standard Brownian motion.
\end{proclaim}
\begin{proclaim}{Theorem D} ({\em Skorokhod embedding theorem})
In a possibly enlarged probability space in which there is a
standard motion $B(t)$, we can redefine the martingale sequence
$\{M_n,\Cal F_n\}$ without changing its distribution and define a
non-decreasing sequence of random variables  $\{\tau_n\}$ such that
$\tau_n$ is $\Cal F_n$ measurable,  $M_n=B(\tau_n)$ and $\ep[\Delta
\tau_m|\Cal F_{m-1}]=\ep[(\Delta M_m)^2|\Cal F_{m-1}]$ a.s.,
$\ep[|\Delta\tau_m|^P|\Cal F_{m-1}]\le C_p\ep[|\Delta M_m|^{2p}|\Cal
F_{m-1}]$ a.s. for $p\ge 1$.
\end{proclaim}

\begin{proclaim}{Theorem E} ({\em Strong approximation})
Let $\{\bm M_n, \Cal F_n; \ge 1\}$ be a martingale sequence in
$\mathscr R^d$ space,  $\bm \Sigma_n =\sum_{m=1}^n
\ep[(\Delta \bm M_m)^{\prime}\Delta \bm M_m |\Cal F_{m-1}]$. Suppose
that there exists constants $0<\epsilon<1$ such that
\begin{equation}\label{eqthE1}\sum_{n=1}^{\infty} \ep \big[\|\Delta\bm M_n\|^2 I\{\|\Delta\bm M_n\|^2\ge
n^{1-\epsilon}\big| \Cal F_{n-1}\big]/n^{1-\epsilon}<\infty\quad
a.s.,
\end{equation}
 and that $\bm T$ is a covariance
matrix which is  measurable with respect to $\Cal F_k$ for some
$k\ge 0$. Then for any $\delta>0$, (possibly in an enlarged
probability space with $\{\bm M_n\}$ being redefined) there exist
$\kappa>0$ and a $d$-dimensional standard Brownin $\bm W(t)$,
independent of $\bm T$, such that
$$ \bm S_n-\bm W(n) \bm T^{1/2}=O(n^{1/2-\kappa})
+O(\alpha_n^{1/2+\delta}) \; a.s.
$$ Here $\alpha_n=\max_{m\le
n}\|\bm \Sigma_m-m\bm T\|$.
\end{proclaim}
If $\sup_n\ep\|\Delta\bm M_n\|^{2+\delta_0}<\infty$ for some
$\delta_0>0$, then the condition (\ref{eqthE1}) is satisfied. The
proofs of Theorem A (a), Theorem C  and Theorem D can be found in
Hall and Heyde (1980). The proof of Theorem A (b) can be found in
Stout (1974). Theorem B can be proved by Theorem D and the LIL of a
Brownian motion. The proof of Theorem E is given in Zhang (2004).
Also, Theorems A-C remain true for martingales in a $\mathscr R^d$
space with some necessary notations changed.

\section{Play-the-winner rule and  Markov chain adaptive design}\label{section3}
\setcounter{equation}{0}

Consider a two-arm clinical trial: two treatments (1 and 2) with
dichotomous response (success and failure).  Patients (subjects) are recruited
into the clinical trial sequentially and respond immediately to
treatments. Zelen (1969) proposed the following design, which is
well known as the play-the-winner (PW) rule: {\em A success on a
particular treatment generates a future trial on the same treatment
with a new patient. A failure on a treatment generates a future
trial on the alternate treatment.} Let $p_i$  be the success
probability of a patient on the treatment $i$, $q_i=1-p_i$, $i=1,2$.
Then
$$ \frac {N_{n,1}}{n}\overset{P}\to \frac{q_2}{q_1+q_2} \;\; \text{
and}\;\; \sqrt{n} (\frac
{N_{n,1}}{n}-\frac{q_2}{q_1+q_2})\overset{\mathscr D}\to
N(0,\sigma_{PW}^2),
$$
 where
$ \sigma_{PW}^2=q_1 q_2(p_1+p_2)/(q_1+q_2)^3$.   Notice that
$q_2/(q_1+q_2)>1/2$ if $p_1>p_2$. So, if treatment 1 is "doing
better", the PW rule favors treatment 1.

Lin, Bai, Chen and Hu (2003) extended the PW rule to a general
Markov Chain adaptive design. Suppose that at the stage $m$, the
treatment 1 is assigned to the $m$th patient. Then the $(m+1)$th
patient will be assigned either treatment 1 or treatment 2 according
certain probabilities, which depend on the response of the $m$th
patient. Let $\alpha_s$ be the probability of assigning the
$(m+1)$th patient to treatment 1, when the response of the $m$th
patient to treatment $1$ is $``$success$"$, and let $\alpha_f$ be
the probability of assigning the $(m+1)$th patient to treatment 1,
when the response of the $m$th patient to treatment $1$ is
$``$failure$"$. Similarly define $\beta_s$ and $\beta_f$ with
treatment 2 instead of treatment 1 in the definitions of $\alpha_s$
and $\alpha_f.$ When $\alpha_s=1,$ $\alpha_f=0,$ $\beta_s=1,$
$\beta_f=0,$ we get Zelen's PW rule. A clinical application
disadvantage of the PW rule is that it is fully deterministic, i.e.,
when the previous results are known, the assignment of the next
subject is fully determined. The Markov chain adaptive design is not
fully deterministic except when the parameter  $\alpha_s,$ $\alpha_f,$
$\beta_s$ and $\beta_f$ take extreme values $0$ and $1$. When taking
$\alpha_s=\alpha_f=\beta_s=\beta_f=1/2$, we get the fully
randomization procedure which allocates patients to each treatment
with a probability $1/2$. The more are the parameters near extreme
values, the more is the procedure being deterministic. When
$\alpha_s<\alpha_f$ and $\beta_s<\beta_f$, the Markov chain adaptive
design is less ethical than the balanced allocation (c.f., Equation
(\ref{eqMCAD1})). The parameters $\alpha_s,$ $\alpha_f,$ $\beta_s$
and $\beta_f$   can be chosen to reflect the trade-off between the
degree of randomness and ethic.

 Let
$p_1(m)=P\{\mbox{success }| X_{m,1}=1\}$ and $p_2(m)=P\{\mbox{success
}| X_{m,1}=0\}$ and $\alpha_m=p_1(m)\alpha_s+(1-p_1(m))\alpha_f,$
$\beta_m=p_2(m)\beta_s+(1-p_2(m))\beta_f.$ Then $\{X_{m,1}\}$ is a
Markov chain with the transition probability matrix
$$\bm P_m=\left(
\begin{array}{ll}
\alpha_m & 1-\alpha_m\\
1-\beta_m & \beta_m
\end{array}
\right).$$ When $p_1(m)=p_1$ and $p_2(m)=p_2$ for all $m$,
$\alpha_n\equiv\alpha,$ $\beta_n\equiv\beta$, and $\{X_{m,1}\}$ is a
homogeneous Markov chain with a stationary distribution $(\mu,
1-\mu)$, where $\mu=(1-\beta)/(2-\alpha-\beta)$. Following from the
central limit theorem for Markov chains we have that
\begin{equation}\label{eqMCAD1}
\frac{N_{n,1}}{n}\overset{P}\to\mu \; \text{ and } \;
\sqrt{n}(\frac{N_{n,1}}{n}-\mu)\overset{\cal D}\to N(0, \sigma^2)
\end{equation}
where
$$
\sigma^2=\Var_{\mu}\{X_{1,1}\}+2\sum_{j=2}^{\infty}\Cov_{\mu}\{X_{1,1},X_{j,1}\}
=\frac{(1-\alpha)(1-\beta)(\alpha+\beta)}{(2-\alpha-\beta)^3}.$$ For
non-homogeneous case,  Lin, Bai, Chen and Hu ( 2003) proved
(\ref{eqMCAD1}) under the condition that
$$
\sum_{m=1}^{\infty}\frac{|p_1(m)-p_1|+|p_2(m)-p_2|}{\sqrt{m}}<\infty.
$$

Lin, Zhang, Cheung and Chan (2005) established the strong
approximation for $N_{n,1}$, from which (\ref{eqMCAD1}) follows
immediately.
\begin{theorem} In a possibly enlarged probability space,  we can redefine the sequence
$\{N_{n,1}\}$ without changing its
distribution, such that
$$
N_{n,1}-n\mu- \sigma W(n)=O((n\log\log n)^{1/4}(\log
n)^{1/2})+O(\Delta_n) \quad a.s., $$
 where $\{W(t)\}$ is a standard
Brownian motion and $\Delta_n=\sum_{m=1}^n (|\alpha_m-\alpha|+|\beta_m-\beta|)$.
\end{theorem}
{\bf Proof.} Write $\Delta m_m=X_{m,1}-\ep[X_{m,1}|\Cal F_{m-1}]$,
where $\Cal F_{m-1}$ is the history sigma field generated by
$X_{m,1},\ldots, X_{m-1,1}$. Then $X_m=\Delta m_m+\ep[X_{m,1}|\Cal
F_{m-1}]=\Delta m_m+1-\beta_m-(1-\alpha_m-\beta_m) X_{m-1,1}$. It
follows that $N_{n,1}=m_n+1-\beta-(1-\alpha-\beta)
N_{n-1,1}+O(\Delta_n)+O(1)$. So,
$$ N_{n,1}=n\mu +m_n/(2-\alpha-\beta)+O(\Delta_n). $$
For the martingale $m_n$, we have $\ep[(\Delta m_m)^2|\Cal F_{m-1}]=(1-\beta_m)\beta_m
+(\alpha_m-\beta_m)(1-\alpha_m-\beta_m)X_{m-1,1}$
and $|\Delta m_m|\le 1$. It follows that
\begin{align*}
 &\sum_{m=1}^n\ep[(\Delta m_m)^2|\Cal F_{m-1}]\\
=&n(1-\beta)\beta +(\alpha -\beta) (1-\alpha -\beta )N_{n-1,1}+O(\Delta_n)\\
 =&n(1-\beta)\beta +n(\alpha -\beta )(1-\alpha -\beta )\mu + \frac{(\alpha -\beta )(1-\alpha -\beta )}{2-\alpha-\beta}m_n+O(\Delta_n)\\
 =&n\frac{(1-\alpha)(1-\beta)(\alpha +\beta )}{2-\alpha-\beta}+O(\Delta_n)+O(\sqrt{n\log\log n}) \\
=:&n\sigma_M^2+O(\Delta_n)+O(\sqrt{n\log\log n}) \;\; a.s.
\end{align*}
due to the LIL (Theorem B). Applying the Skorokhod embedding theorem (Theorem D), we can write $m_m=B(\tau_m)$. Notice
$\sum_{m=1}^n\big(\Delta \tau_m-\ep[\Delta \tau_m|\Cal F_{m-1}])$ is also a martingale sequence. By the LIL, we conclude that
\begin{align*}
\tau_n=&\sum_{m=1}^n \ep[\Delta \tau_m|\Cal F_{m-1}]+O(\sqrt{n\log\log n})\\
=&\sum_{m=1}^n \ep[(\Delta m_m)^2|\Cal F_{m-1}]+O(\sqrt{n\log\log n})\\
=&n\sigma_M^2+O(\Delta_n)+O(\sqrt{n\log\log n}) \;\; a.s.
\end{align*}
It follows that
\begin{align*} m_n=&B(\tau_n)=B(n\sigma_M^2)+O\Big(\sqrt{(\Delta_n+\sqrt{n\log\log n})\log n}\Big)\\
=&B(n\sigma_M^2)+O(\Delta_n)+O((n\log\log n)^{1/4}(\log
n)^{1/2}) \;\; a.s.
\end{align*}
by the sample properties of a Brownian motion (c.f., Cs\" org\H o and R\'ev\'sz (1980)).
The proof is now completed by letting $W(t)=B(t\sigma_M^2)/\sigma_M$ and noticing that $\sigma_M^2/(2-\alpha-\beta)^2=\sigma^2$.
$\Box$.

\smallskip
 For the multi-treatment
case, we let
$p_i(m)=\pr(\text{success }| X_{m,i}=1\}$. Assume the transition
probability matrix of the Markov chain $\{\bm X_n\}$ is $\bm
H_n=\{H_{ij}(n)\}$ which is a function of $p_i(n)$, $i=1,\ldots,K$,
i.e, $\ep[\bm X_{n+1} | \bm X_n]= \bm X_n \bm H_n$. By using Theorem E instead of the Skorokhod embedding theorem,
Zhang (2004) showed that
$\bm N_n$ can be approximated by a multi-dimensional Browian motion:
$$
\bm N_n-n\bm v- \bm W(n) \bm\Sigma
=o(n^{1/2-\kappa})+O\big(\sum_{k=1}^n\|\bm H_k-\bm H\|\big)
 \quad a.s.,
$$
where $\kappa>0$, $\bm v=(v_1, \cdots, v_K)$ is the left eigenvector
corresponding to the largest eigenvalue of $\bm H$ with
$v_1+\ldots+v_K=1$, $\{\bm W(t)\}$ is a $K$-dimensional standard
Brownian motion. In particular, we have asymptotic normality, if
\begin{equation}\label{co3}
\sum_{k=1}^n\|\bm H_k-\bm H\|=o(n^{1/2}).
\end{equation}
Zhang (2006) studied a kind of non-humongous Markov chain designs, in which  $\bm H_n$
is a function of an  estimated unknown parameters $\widehat{\bm \theta}$.
In this case the condition (\ref{co3}) is not satisfied.

\section{Randomized play-the-winner rule and  urn
models}\label{section4} \setcounter{equation}{0}

  To overcome the drawback that the PW rule is fully
deterministic, Wei and Durham (1978) introduced the following
randomized play-the-winner (RPW) rule: {\em We start with
$(\alpha,\beta)$ balls (type $1$ and $2$ respectively) in the urn.
If a type $k$ ball is drawn, a patient is assigned to the treatment
$k$, $k=1,2$. The ball is replaced and the patient's response is
observed. A success on the treatment $1$ or a failure on the
treatment $2$ generates a type $1$ ball in the urn; A success on the
treatment $2$ or a failure on the treatment $1$ generates a type $2$
ball in the urn.} Let $Y_{n,1}$ ($Y_{n,2}$) be number of balls of
type 1 (2) after $n$ stage.
 From the results of Athreya and Karlin (1968), we have
$$ \frac {Y_{n,1}}{Y_{n,1}+Y_{n,2}}\to \frac{q_2}{q_1+q_2}\; a.s.
\quad\text{ and } \quad
 \frac {N_{n,1}}{n}\to \frac{q_2}{q_1+q_2} \; a.s.. $$
The limiting proportion is the same as that of the PW rule. We refer to it as urn proportion.
When $p_1+p_2<1.5$ (or $q_1+q_2>0.5$), we have the following
asymptotic normality:
$$
\sqrt{n} (\frac {Y_{n,1}}{Y_{n,1}+Y_{n,2}}-\frac{q_2}{q_1+q_2})\overset{\Cal D}\to
N\big(0,\frac{q_1 q_2}{(2(q_1+q_2)-1)(q_1+q_2)^2}\big)
$$
and $ \sqrt{n} (N_{n,1}/n-q_2/(q_1+q_2))\overset{\Cal D}\to
N(0,\sigma_{RPW}^2),
$
where
\begin{equation}\label{eqvarofRPW}
 \sigma_{RPW}^2=\frac{q_1 q_2[5-2(q_1+q_2)]}{[2(q_1+q_2)-1](q_1+q_2)^2}.
\end{equation}
 The asymptotic normality was first given in Smythe and Rosenberger
(1995). When $q_1+q_2<0.5$, the limiting distributions of both the
urn composition and the allocation proportion  are unknown. The RPW
rule has the same limiting allocation proportion as the PW rule. But
the asymptotic variability is much larger.

\smallskip
 As multi-treatment extensions of the RPW rule, one large
family of randomized adaptive designs can be developed from the
generalized Polya urn (GPU) model. Urn models have also long been
recognized as valuable mathematical apparatus in many areas
including physical science, biological science, engineering,
information science, the study of economic behaviors,  etc.

Consider an urn containing balls of K types. Initially, the urn
contains ${\bm Y}_0=(Y_{0,1}, \cdots, Y_{0,K})$ balls, where
$Y_{0,k}$ denotes the number of balls of type $k$, $k=1,\cdots,K$. A
ball is drawn at random. Its type is observed and the ball is then
replaced. At the $m$th stage, following a type $k$ drawn, $D_{kj}(m)
$ ($\ge 0$) balls of type $j$, for $j=1,\cdots,K$, are added to the
urn. $D_{ij}(m)$ is  a random function of the response $\xi_{m,k}$
of the $m$-th subject on treatment $k$. The expectation of the total
numbers of balls added in each stage is assumed to be the same (say
$\gamma$), so
$$\sum_{j=1}^K\ep\{ D_{kj}(m)|\Cal F_{m-1}\}=\gamma,\quad  k,j=1, \cdots, K,
m=1,2,\cdots, $$ where $\Cal F_{m-1}$ is the history sigma field.
Without loss generality, we can assume $\gamma=1$.
 Let ${\bm H_m}$ be the matrix comprising element
$\{h_{ij}(m)= \ep[D_{kj}(m)|\Cal F_{m-1}]\}$ and $\bm
D_m=\{D_{kj}(m)\}$. We refer to ${\bm D_m}$ as the adding
 rules and ${\bm H_m}$ as the  design matrices.
 If $\bm D_m$, $m=1,2,\ldots,$  are i.i.d., then $\bm H_m =\bm
H$ for all $n$. In this case  the model is said to be homogenous. In
general, it is assumed that $\bm H_m \to \bm H$. Let ${\bm
Y}_n=(Y_{n,1}, \cdots, Y_{n,K})$, where $Y_{n,k}$ represents the
number of balls in the urn of type $k$ after $n$th stage. And let
$\bm v=(v_1, \cdots, v_K)$ be the left eigenvector corresponding to
the largest eigenvalue of $\bm H$ with $v_1+\ldots+v_K=1$. Then
$v_k$ is just the limiting proportion of both the patients assigned
to treatment $k$ and the type $k$ balls in the urn. Bai, Hu and
Zhang (2002), Hu and Zhang (2001) obtained the asymptotic properties
via the strong approximation.
\begin{theorem}\label{urnth} Suppose that $\{\bm D_m\}$ is sequence of i.i.d. random matrices
with $\sup_m\ep\|\bm D_m\|^{2+\delta}<\infty$. Let
$\lambda_1=\gamma=1, \lambda_2,\ldots, \lambda_K$ be the eigenvalues
of $\bm H$, and
$\lambda=\max\{Re(\lambda_2),\ldots,Re(\lambda_K)\}$. If
$\lambda<1$, then
\begin{equation}\label{equrncon}
\frac{N_{n,i}}{n}\to v_i\quad a.s.
\; \text{ and } \;
\frac{Y_{n,i}}{\sum_{j=1}^{K} Y_{nj}}\to v_i\quad a.s.
\end{equation}
If $\lambda<1/2$, then
\begin{equation}\label{urnnormal}
\sqrt{n}(\frac{\bm Y_n}{n}-\bm v)\overset{\Cal D}\to
N(0,\bm\Sigma)\;\; \text{ and }\;\;
\sqrt{n}(\frac{\bm N_n}{n}-\bm v)\overset{\Cal D}\to
N(0,\bm\Sigma^{\ast}).
\end{equation}
\end{theorem}
{\bf Proof.} Write $\bk 1=(1,\ldots, 1)$, $|\bm Y_n|=\sum_{k=1}^K
Y_{nk}$, $\widetilde{\bk H}=\bm H-\bm 1^{\prime}\bm v$, $\bm
M_n=\sum_{m=1}^N \bm X_m (\bm D_m-\bm H)$, $\bm m_m=\sum_{m=1}^n
(\bm X_m-\ep[\bm X_m|\Cal F_{m-1}])$. Then $\bm H\bm 1^{\prime}=\bm
1^{\prime}$, $\bm m_n\bm 1^{\prime}=\bm 0$, and the eigenvalues of
$\widetilde{\bk H}$ are $0,\lambda_2,\ldots, \lambda_K$. We have the
following lemma on matrices, the proof of which can be founded in Hu
and Zhang (2004a).
\begin{lemma}\label{urnlem} If $\Delta \bm Q_n=\Delta \bm P_n+\bm Q_{n-1}\widetilde{\bm H}/(n-1)$, $n\ge 2$, then
$$ \|\bm Q_n\|=O(\|\bm P_n\|)+\sum_{m=1}^n\frac{O(\|\bm P_m\|) }{m}\big(n/m)^{\lambda}\log^{\nu-1}(n/m), $$
where $\nu$ is the degree of the second largest eigenvalue of $\bm H$.
\end{lemma}
We prove (\ref{urnnormal}) only. Notice
\begin{align}\label{equrniterative}
&\Delta (\bm Y_n-n\bm v)=\bm X_n\bm D_n-\bm v
 =\Delta \bm M_n +\Delta \bm m_n\bm H+\frac{\bm Y_{n-1}}{|\bm Y_{n-1}|}\bm H -\bm v
 \nonumber\\
&=\frac{\bm Y_{n-1}-(n-1)\bm v}{n-1}\widetilde{\bm H}+\Delta \bm M_n
+\Delta \bm m_n\widetilde{\bm H}
\nonumber\\
&\qquad +\left(1-\frac{|\bm Y_{n-1}|}{ n-1 }\right)\left(\frac{\bm
Y_{n-1}}{|\bm Y_{n-1}|}-\bm v\right)\widetilde{\bm H}.
\end{align}
Multiplying $\bm 1^{\prime}$ yields
$ |\bm Y_n|-n=\Delta \bm M_n\bm 1^{\prime}+|\bm Y_0|=O(\sqrt{n\log\log n})\;\; a.s. $
due to the LIL (Theorem B). It follows that
$$  \bm M_n + \bm m_n\bm H
+\sum_{m=1}^{n-1}\left(1-\frac{|\bm Y_m|}{m }\right)\left(\frac{\bm
Y_m}{|\bm Y_m|}-\bm v\right)\widetilde{\bm H} =O(\sqrt{n\log\log
n})\;a.s. $$ By applying Lemma \ref{urnlem} and noticing $\lambda
<1/2$, we obtain $\bm Y_n-n\bm v=O(\sqrt{n\log\log n})$\;a.s. Write
$\bm D_1^{(k)}=(D_{k1},\ldots, D_{kK})$, $\bm \Sigma_1=diag(\bm v)-\bm
v^{\prime}\bm v$ and  $\bm \Sigma_2=\sum_{k=1}^K
v_k \Var\{\bm D_1^{(k)}\}$. For the martingale $(\bm M_n, \bm m_n)$, we have
$\sup_n\ep\|\Delta\bm M_n\|^{2+\delta}<\infty$, $\|\Delta \bm
m_n\|\le K$, $\ep[(\Delta\bm M_n)^{\prime}\Delta\bm m_n|\Cal
F_{m-1}]=\bm 0$  and
\begin{align*}
\sum_{m=1}^n\ep\Big[(\Delta\bm m_m)^{\prime}\Delta\bm m_m|\Cal
F_{m-1}\Big] =&\sum_{m=1}^n[diag\big(\frac{\bm Y_m}{|\bm Y_m|}\big)
-\big(\frac{\bm Y_m}{|\bm Y_m|}\big)^{\prime}\frac{\bm Y_m}{|\bm Y_m|}]\\
=&n\bm \Sigma_1+O(\sqrt{n\log\log n})\;\; a.s.,\\
\sum_{m=1}^n\ep\Big[(\Delta\bm M_m)^{\prime}\Delta\bm M_m|\Cal
F_{m-1}\Big]
=&\sum_{m=1}^n\sum_{k=1}^K\Var\{\bm D_1^{(k)}\}Y_{mk}/|\bm Y_m|\\
=&n\bm \Sigma_2+O(\sqrt{n\log\log n})\;\; a.s.
\end{align*}
By the strong approximation (Theorem E), there are two independent $d$-dimensional standard Brownian motions $\bm W_1(t)$ and
$\bm W_2(t)$ such that for some $\kappa>0$,
\begin{equation}\label{equrnappofM}
\bm m_n=\bm W_1(n)\bm\Sigma_1^{1/2}+o(n^{1/2-\kappa}) \; a.s., \;
 \bm M_n=\bm W_2(n)\bm\Sigma_2^{1/2}+o(n^{1/2-\kappa}) \; a.s.
\end{equation}
Without loss generality, we assume $1/2-\kappa>\lambda$ and $1/2-\kappa>1/4$.
Let $\bm G_i(t)$ be the solution of the equation
\begin{equation}\label{equrnsolution}
 \bm G_i(t)=\int_0^t \frac{\bm G_i(x)}{x}\widetilde{\bm H} dx +\bm W_i(t)\bm\Sigma_i^{1/2},\;\; \bm G_i(0)=0,
 \end{equation}
$i=1,2$. Notice $\{[\bm W_i(T(\cdot+s))-\bm
W_i(T\cdot)]/\sqrt{T}\}\overset{\mathscr D}= \{\bm W_i(\cdot)\}$. It
can be checked that $\bm G_i(t)$ is a Gaussian process with
stationary increments and $\Var\{\bm G_i(t)\}=t\Var\{\bm G_i(1)\}$,
and
\begin{equation}\label{equrnintofG}
 \int_0^n \frac{\bm G_i(x)}{x}  dx=\sum_{m=1}^{n-1}\frac{\bm G_i(m)}{m} +O((\log n)^{3/2})\;\; a.s.
 \end{equation}
Combing (\ref{equrniterative})--(\ref{equrnintofG}) yields
$$ \bm Y_n-n\bm v-\bm G_2(n) -\bm G_1(n) \widetilde{\bm H}
=\sum_{m=1}^{n-1}\frac{\bm Y_m-m\bm v-\bm G_2(m) -\bm G_1(m) \widetilde{\bm H}}{m}
\widetilde{\bm H}+\bm P_n,
$$
where
\begin{align*}
&\bm P_n=O((\log n)^{3/2})+o(n^{1/2-\kappa})
+\sum_{m=1}^{n-1}\left(1-\frac{|\bm Y_m|}{m }\right)\left(\frac{\bm Y_m}{|\bm Y_m|}-\bm v\right)\widetilde{\bm H}\\
=&O((\log n)^{3/2})+o(n^{1/2-\kappa})+\sum_{m=1}^{n-1}\left(\sqrt{(\log\log m)/m}\right)^2
=o(n^{1/2-\kappa})\;\; a.s.
\end{align*}
Now, by applying Lemma \ref{urnlem} we conclude that
$$\bm Y_n-n\bm v=\bm G_2(n) +\bm G_1(n) \widetilde{\bm H}+
o(n^{1/2-\kappa})\;\; a.s. $$
Finally,
\begin{align*}
\bm N_n- n\bm v=& \bm m_n+\sum_{m=0}^{n-1}\frac{\bm Y_m}{|\bm
Y_m|}-n\bm v =
\bm m_n+\sum_{m=0}^{n-1}\left(\frac{\bm Y_m}{|\bm Y_m|}-\bm v\right)(\bm I-\bm 1^{\prime}\bm v)\\
=&
\bm m_n+\sum_{m=1}^{n-1}\frac{\bm Y_m-m\bm v}{m}(\bm I-\bm 1^{\prime}\bm v)\\
&+\sum_{m=1}^{n-1}\left(1-\frac{|\bm Y_m|}{m}\right)\left(\frac{\bm Y_m}{|\bm Y_m|}-\bm v\right)(\bm I-\bm 1^{\prime}\bm v)
+O(1)\\
=&
\bm W_1(n)+\sum_{m=1}^{n-1}\frac{\bm G_2(m) +\bm G_1(m) \widetilde{\bm H}}{m}(\bm I-\bm 1^{\prime}\bm v)
+o(n^{1/2-\kappa})\\
=&
\bm W_1(n)+\int_0^n\frac{\bm G_2(x) +\bm G_1(x) \widetilde{\bm H}}{x}dx(\bm I-\bm 1^{\prime}\bm v)
+o(n^{1/2-\kappa})\\
=&
\bm G_1(n)+\int_0^n\frac{\bm G_2(x)}{x}dx(\bm I-\bm 1^{\prime}\bm v)
+o(n^{1/2-\kappa})\;\; a.s.
\end{align*}
(\ref{urnnormal}) is now proved, where $\bm\Sigma$ and $\bm \Sigma^{\ast}$ are the variance-covariance matrices
 of normal random variables $\bm G_2(1)+\bm G_1(1)\widetilde{\bm H}$ and
 $\bm G_1(1)+\int_0^1 \bm G_2(x)/x\, dx$ $ (\bm I - \bm 1^{\prime}\bm v)$, respectively. For the details of
 specifying  and estimating the variance-covariance matrices, one can refer to Bai and Hu (2005), Hu and Zhang (2004b),
 Zhang, Hu and Cheung  (2006). $\Box$

 The interested one can check that for the RPW rule, the equation (\ref{equrnsolution}) reduces to
 $$ G_{i,1}(t)=\lambda \int_0^t \frac{G_{i,1}(x)}{x}dx +\sigma_i W_i(t),\; G_{i,1}(0)=0, $$
 $G_{i,2}(t)=-G_{i,1}(t)$,  and the solution is
 $ G_{i,1}=\sigma_it^{\lambda}\int_0^t x^{-\lambda} dW_i(x)$, $i=1,2$, where $\lambda=1-q_1-q_2$, $\sigma_1^2
 =v_1v_2=q_1q_2/(q_1+q_2)^2$ and $\sigma_2^2=v_1p_1q_1+v_2p_2q_2=q_1q_2(p_1+p_2)/(q_1+q_2)$.

\smallskip
When only the second moment of $\bm D_m$ is assumed to be finite, we
can show (\ref{urnnormal}) with a similar argument by applying the
weak convergence of martingales (Theorem C) instead of the strong
approximation. For details one can refer to Hu and Zhang (2001).
Janson (2004) studied the urn models by embedding them to continuous
branching processes and established the asymptotic normality in a
different way.

\begin{example} \label{exampleWei79} As a multi-treatment extension of the RPW rule, Wei
(1979) proposed a GUP to allocate subjects, in which the urn is
updated in  the following way: at the $n$th stage, if a subject is
assigned to treatment $k$ and cured, then a type $k$ ball is added
to the urn, otherwise, if treatment $k$ for a subject fails, then
$\frac 1{K-1}$ balls are added to the urn for each of the other
$K-1$ treatments. In this urn model, $\bm H =\{h_{kj},
k,j=1,\ldots,k\}$, where $h_{kk} =p_k$ and $h_{kj} =q_k/(K-1)$
($j\ne k$), and $p_k$ is the successful probability of treatment
$k$, $q_k=1-p_k$. So $ v_k=(1/q_K)/\sum_{j=1}^{K}(1/q_j)$.
\end{example}

For the non-homogenous case,
 Bai and Hu (2005) obtained (\ref{urnnormal}) under the condition that
\begin{equation}\label{conditionurn}
 \sum_{m=1}^{\infty}\frac{\|\bm H_m-\bm H\|}{\sqrt m}<\infty \;a.s.
\end{equation}
This condition can be weakened to $\sum_{m=1}^n \|\bm H_m-\bm
H\|=o(\sqrt{n})$ a.s. by use the argument in the above proof. An
applicable class of non-homogenous urn models is the sequential
estimation-adjusted urn (SEU) model, in which the urn is updated
according to the current response and the current estimate of an
unknown parameter, and so $\bm H_m=\bm H(\widehat{\bm\theta}_m)$ is
a function of the estimator. In this case, the fastest convergence
rate of $\bm H_m$ is $O_P(\sqrt{m})$ and the condition
(\ref{conditionurn}) is not satisfied. Zhang, Hu and Cheung  (2006)
established the asymptotic properties of SEU models.
\begin{example} Bai, Hu and Shen (2002) proposed a GUP to allocate subjects,
in which the urn is updated in  the following way: at the $m$th
stage, if a subject is assigned to treatment $k$ and cured, then a
type $k$ ball is added to the urn. If a failure, then
$\frac{\widehat{p}_{m-1,j}}{\sum_{i\ne k} \widehat{p}_{m-1,i}}$
balls of each type $j\ne k$ are added. Where
$\widehat{p}_{m-1,j}=(S_{m-1,i}+1)/(N_{m-1,j}+1)$, and $S_{m-1,j}$
is the number of successes of treatment $j$ in previous $m-1$
stages. This model is a SEU model with $\bm H(\bm x)=\{h_{ij}(\bm
x);i,j=1,\ldots, K\}$, where $h_{kk}(\bm x)=p_k$ and $h_{kj}(\bm
x)=q_k x_j/\sum_{i\ne k}x_i$ ($j\ne k$).
\end{example}
More examples and applications of SEU models can be found in Zhang,
Hu and Chueng (2006), in which how to defined a SEU model by using
the information of distribution parameters to target a pre-specified
limiting allocation proportion is discussed in details.

\section{Drop-the-loser rule}\label{sectionad4}
\setcounter{equation}{0}

The asymptotic normality for the urn models can be obtained only
when the condition $\lambda\le 1/2$ is satisfied. This is a very
strict condition. Even in the case of $K=3$, it is hard to be
satisfied and to check it is not a easy work. Also, when $\lambda$
is close  or exceeds $1/2$, the variability of an urn model is
extremely high. Ivanova (2003) proposed a drop-the-loser (DL) rule
which has the same limiting proportion as Wei (1979)'s rule (See
Example \ref{exampleWei79}) but has much smaller variability.
Consider  an urn containing balls of $K+1$ types, type $0$,
$1,\ldots K$, when comparing $K$ treatments. A ball is drawn at
random. If it is type $k$, $k=1,\ldots,K$, the corresponding
treatment is assigned and the subject's response is observed. If the
response is a success, the ball is replaced and the urn remains
unchanged. If   a failure, the ball is not replaced. When a type $0$
ball is drawn, no subject is treated, and the ball is return to the
urn together with one ball of   each type $k$, $k=1,\ldots, K$.
Ivanova (2003, 2006) established the asymptotic normality after
embedding the urn process to a death-and-immigration process. Here
we give the strong approximation.
\begin{theorem} There is a $K$-dimensional standard Brownian motion
$\bm W(t)$ such that
$$ \bm N_n-n\bm v=\bm W(n)
diag\Big(\sqrt{\frac{v_1p_1}{q_1}},\ldots,\sqrt{\frac{v_Kp_K}{q_K}}\Big)(\bm
I-\bm 1^{\prime}\bm v)+o(n^{1/2-\kappa}) \;\; a.s., $$ for some
$\kappa>0$, where $v_k=(1/q_k)/\sum_{j=1}^k(1/q_j)$, $k=1,\ldots,K$.
Hence
$$\sqrt{n}(\bm N_n/n-\bk v)\overset{\mathscr D}\to N\Big(\bm 0, (\bm I-\bm v\bm 1^{\prime})
diag\big( \frac{v_1p_1}{q_1} ,\ldots, \frac{v_Kp_K}{q_K}\big)(\bm I-\bm 1^{\prime}\bm v)\Big).$$
\end{theorem}
In particular, in the two-treatment case, $\sqrt{n}(\bm N_n/n-\bk v)\overset{\mathscr D}\to N(0,\sigma_{DL}^2)$
with $\sigma_{DL}^2=q_1q_2(p_1+p_2)/(q_1+q_2)^3$, the same as the $\sigma_{PW}^2$.
For the generalizations of the DL rule and their applications,
one can refer to Zhang, Chan, Cheung and Hu (2007), Sun, Cheung and Zhang (2007).

\smallskip
{\bf Proof of the theorem.} Let $\bm Z_m=(Z_{m,0},\ldots, Z_{m,K})$ be the urn compositions after the $m$-the assignment.
And let $\mu_m$ be the number of draws of type 0 balls between the $(m-1)$-th assigment and the $m$-th  assignment. Remember
that when a  type 0 ball is drawn, we add one ball of each treatment type, and when a treatment type ball is drawn, it is replaced
only when the response is a success. So
$$ Z_{m,k}-Z_{m-1,k}=\mu_m+X_{m,k}(\xi_{m,k}-1)=\mu_m-X_{m,k}q_k+X_{m,k}(\xi_{m,k}-p_k), $$
where $\xi_{m,k}=1$ if the response of the $m$-th subject on
treatment $k$ is a success, and $0$ if failure. Let $
M_{n,k}=\sum_{m=1}^nX_{m,k}(\xi_{m,k}-p_k)$, $\bm
M_n=(M_{n,1},\ldots, M_{n,K})$,  and $\Cal A_m$ be the sigma field
generated by $\xi_{1,k},\cdots,\xi_{m,k}$, $k=1,\ldots,K$, and $\bm
X_1,\ldots, \bm X_m, \bm X_{m+1}$. Then $\{\bm M_n,\Cal A_m\}$ is a
martingale. It follows that $ Z_{n,k}-Z_{0,0}=\sum_{m=1}^n
\mu_m-N_{n,k}q_k+M_{m,k}$. We can prove that
$Z_{n,k}=o(n^{1/2-\delta_0})$ a.s. for some $\delta_0>0$. For
details of the proof, we refer to Zhang, Hu, Chueng and Chan
(2006b), Sun, Cheung and Zhang (2007). Hence
$$ N_{n,k}=\sum_{m=1}^n \mu_m/q_k+M_{m,k}/q_k+o(n^{1/2-\delta_0})\;\; a.s., \; k=1,\ldots, K, $$
which, together with the fact $N_{n,1}+\cdots+N_{n,K}=n$, implies
$$ \bm N_n-n\bm v=\bm M_n diag(1/q_1,\ldots,1/q_K)(\bm I-\bm 1^{\prime}\bm v)+o(n^{1/2-\delta_0})\;\; a.s.$$
So, $\bm N_n-n\bm v=O(\sqrt{n\log\log n})$ a.s. by the LIL (Theorem B).
On the other hand, for the martingale $\{\bm M_n\}$ we have
\begin{align*}
 &\sum_{m=1}^n \ep[(\Delta\bm M_n)^{\prime}\Delta\bm M_n|\Cal F_{m-1}]=
 diag(N_{n,1}p_1q_1,\ldots, N_{n,K}p_kq_k)\\
&\quad  =n\; diag(v_1p_1q_1,\ldots, v_Kp_kq_k)+O(\sqrt{n\log\log n})\; a.s.
 \end{align*}
By applying the strong approximation (Theorem E), we can define a
$K$-dimensional Brownian motion $\bm W(t)$ such that
$$ \bm M_n=\bm W(n)diag(\sqrt{v_1p_1q_1},\ldots, \sqrt{v_Kp_kq_k})+o(n^{1/2-\kappa})\;\; a.s. $$
The proof is now completed. $\Box$

\section{Doubly adaptive biased coin  designs}\label{section5}

\setcounter{equation}{0}

The PW rule and urn model designs are a kind of design-driven adaptive designs, which  are
  constructed with intuitive motivation.
  However, clinical trials are usual complex experiments on humans with multiple, often competing, objectives,
  including maximizing power to detect clinically relevant differences in treatment outcomes,
  maximizing the individual patient's personal experience while treated in the trial, and minimizing the
  total monetary cost of trial.  These and other objectives can be defined in terms
  of optimization of function of the trial's parameters, the optimal allocation proportion is often
a function of unknown parameters. Take a binary response clinical
trail with two treatments $1$ and $2$ as an example. The well known
Neyman proportion is
$$ \rho(p_1,p_2)=: \frac{n_1}{n_1+n_2} =\frac{\sqrt{p_1q_1}}{\sqrt{p_1q_1}+\sqrt{p_2q_2}}, $$
where $p_k$ ($q_k$) is the probability of success (failure) of a trial
treatment $k$, $n_k$ is the number of subjects assigned to treatment
$k$,  $k=1,2$, The Neyman proportion maximizes the power of a test
of the simple difference $p_1-p_2$ for fixed sample size $n$. But if
we implement Neyman allocation, when $p_1+p_2>1$,
 we will  assign more subjects to the inferior treatment, which will compromise the ethical objective.
 Rosenberger, et al (2001) discussed another important optimization criteria that
 minimize the expected number of treatment  failures, $n_1 q_1+n_2 q_2$, for fixed the variance,
 $p_1q_1/n_1+p_2q_2/n_2$, of the   statistic
 $\widehat{p}_1-\widehat{p}_2$ under an alternative hypothesis $p_1\ne p_2$.  This leads the
 optimal
 proportion as follows.
$$ \rho(p_1,p_2)= \frac{\sqrt{p_1 }}{\sqrt{p_1 }+\sqrt{p_2 }}. $$
Many simulation studies have validated that a adaptive design with
this proportion as its target performs very  satisfactorily for both
the consideration of ethic and the test of power. For other
optimization criteria one can refer to Jennison and Turnbull (2000)
and Rosenberger, et al. (2001). Both Neyman allocation and the
allocation of Rosenberger,et al (2001) cannot be implemented
directly in a clinical trail, because we do not know the values of
$p_1$ and $p_2$. In this section, we introduce an important  class
of adaptive designs, doubly adaptive biased coin designs (DBCD),
which is first proposed by Eisele (1994) and Eisele and Woodroofe
(1995),  to target a pre-specified allocation proportion.

Consider a clinical with $K$ treatments. The outcome of a subject on
treatment $k$ has a distribution $f_k(\cdot|\theta_k)$. Write $\bm
\theta=(\theta_1,\ldots,\theta_K)$. The pre-specified allocation
proportion is $\bm v=\bm \rho(\bm \theta)=(\rho_1(\bm
\theta),\ldots,\rho_m(\bm\theta))$. Here $\bm\rho(\bm y)$ is assumed
to be continuous function on the parameter space, taking the values
on $(0,1)^{\otimes k}$ and twice differentiable at the true value of
the parameter $\bm\theta$. A multi-treatment DBCD proposed by Hu and
Zhang (2004a) is defined as follows.

 To start, allocate $M$ subjects to each treatment. At stage $m$, suppose $m-1$($\ge MK$)
 subjects are allocated and the outcomes, $N_{m-1,k}$ outcomes of treatment $k$, $k=1,\ldots, K$, are observed.
 Let $\widehat{\theta}_{m-1,k}$ be the MLE of the parameter $\theta_k$, $k=1,\ldots,
 K$.
Write $\widehat{\bm
\theta}_{m-1}=(\widehat{\theta}_{m-1,1},\cdots,\widehat{\theta}_{m-1,K})$,
and let $\widehat{\bm \rho}_{m-1}=\rho(\widehat{\bm \theta}_{m-1})$
be the current estimate of the target allocation proportion. Now,
the $m$-th subject is allocated to treatment $k$ with a probability:
\begin{equation}\label{eqallocationDBCD}
P_{m,k}=:\pr(X_{m,k}=1|\Cal F_{m-1})=g_k\big(\frac{\bm
N_{m-1}}{m-1}, \widehat{\rho}_{m-1}\big),
\end{equation}
$k=1,\ldots, K$. Here $\bm g(\bm x,\bm y)=(g_1(\bm x,\bm y),\ldots, g_K(\bm x,\bm y)): (0,1)^{\otimes 2k}\to
(0,1)^{\otimes k}$ is the allocation function. Write $\bm P_m=(P_{m,1},\ldots, P_{m,K})$.

\begin{theorem}\label{thDBCD}Suppose the distributions $f_1(\cdot|\theta_1),\ldots,f_K(\cdot|\theta_K)$ follow an exponential family.
Let $\bk g(\bm x,\bm y)$ be defined as
\begin{equation}\label{eqallocationfunction} g_k(\bm x,\bm y)=\frac{y_k\left(\frac{y_k}{x_k}\right)^{\gamma}}{\sum_{j=1}^Ky_j\left(\frac{y_j}{x_j}\right)^\gamma},
\;\; k=1,\ldots, K; \quad \gamma\ge 0.
\end{equation}
Then
$$ \bm N_n-n\bm v= O(\sqrt{n\log\log n})\; a.s.\;\; \text{ and }
\;\; \sqrt{n}(\bm N_n/n-\bm v)\overset{\mathscr D}\to N(\bm 0,\bm \Sigma),$$
where
$$ \bm \Sigma=\bm \Sigma_{\bm \rho}+\frac{1}{1+2\gamma} (diag(\bm v)-\bm v^{\prime}\bm v+\bm \Sigma_{\bm\rho}), $$
$$\bm \Sigma_{\bm \rho}=\left(\frac{\partial \bm\rho}{\partial \bm\theta}\right)^{\prime}
diag\left((v_1 I_1(\theta_1))^{-1},\ldots,v_K I_1(\theta_K))^{-1}\right)\frac{\partial \bm\rho}{\partial \bm\theta}
$$
and $I_k(\theta_k)$ is the Fisher information function for a single
observation on treatment $k$.
\end{theorem}
In particular, for the two-treatment case, suppose the targeted
allocation proportion of treatment 1 is $
\rho=\rho(\theta_1,\theta_2)$. Then
$ \sqrt{n}\big(N_{n,1}/n-\rho
\big)\overset{\mathscr D}\to N(0,\sigma^2_{DBCD})$
where
$$\sigma^2_{DBCD}=\sigma^2_{\rho}+\frac{1}{1+2\gamma} \big\{\rho (1-\rho )
+\sigma^2_{\rho} \big\}$$
and $\sigma_{\rho}^2=\big(I_1(\theta_1)\rho\big)^{-1}(\partial\rho/\partial \theta_1)^2
+\big(I_2(\theta_2)(1-\rho)\big)^{-1}(\partial\rho/\partial \theta_2)^2 $.

\begin{example} Consider the binary response clinical
trail with two treatments. For the urn proportion
$\rho=q_2/(q_1+q_2)$,
 $$ \sigma_{DBCD}^2=\frac{q_1q_2(p_1+p_2)}{(q_1+q_2)^3}
 +\frac{2q_1q_2}{(1+2\gamma)(q_1+q_2)^3}. $$
For the Neyman proportion $\rho=\sqrt{p_1q_1}/(\sqrt{p_1q_1}+\sqrt{p_2q_2})$, $\sigma_{DBCD}^2 =$
\begin{eqnarray*}&\frac{\sqrt{p_1q_1p_2q_2}}{(1+2\gamma)(\sqrt{p_1q_1}+\sqrt{p_2q_2})^2}+&\\
&\frac{1+\gamma}{2(1+2\gamma)(\sqrt{p_1q_1}+\sqrt{p_2q_2})^3}
\left(\frac{p_2q_2(q_1-p_1)^2}{\sqrt{p_1q_1}}
+\frac{p_1q_1(q_2-p_2)^2}{\sqrt{p_2q_2}}\right).&
\end{eqnarray*}
 For Rosenberger, et al's
proportion $\rho=\sqrt{p_1}/(\sqrt{p_1}+\sqrt{p_2})$, $\sigma_{DBCD}^2 =$
$$\frac{\sqrt{p_1p_2}}{(1+2\gamma)(\sqrt{p_1}+\sqrt{p_2})^2}+
\frac{1+\gamma}{2(1+2\gamma)(\sqrt{p_1}+\sqrt{p_2})^3}\left(\frac{p_2q_1}{\sqrt{p_1}}+\frac{p_1q_2}{\sqrt{p_2}}\right).
$$
\end{example}

From Theorem \ref{thDBCD}, we find that the asymptotic variability
is a decreasing function of parameter $\gamma$. However, the degree
of randomness of the design  decreases when $\gamma$ increases,
because, as the value of $\gamma$ becomes larger, the
allocation probabilities shift faster to extreme values $0$ and $1$   if
there is a bias between the current sample allocation  and the
estimated target. When $\gamma=\infty$, the variability of the
procedure is minimized, but the procedure is completely predictable.
The parameter $\gamma$ can be chosen to reflect the trade-off
between the degree of randomness and the variability.

The allocation function defined in (\ref{eqallocationfunction}) is
very special though it has fine properties.
 For results for general allocation function $\bm g(\cdot,\cdot)$, one can refer to Hu and Zhang
 (2004a).

\smallskip
{\bf Proof the Theorem.} Let $\xi_{m,k}$, $m=1,2,\ldots$, be i.i.d.
random variables, which represent  the outcomes on treatment $k$,
$k=1,\ldots,K$. In clinical trial, only $X_{m,k}\xi_{m,k}$s are
observed. Write $\bm \xi_m=(\xi_{m,1},\ldots,\xi_{m,k})$. For
simplifying the proof, we assume that $\theta_k=\ep\xi_{m,k}$ is the
mean of the outcomes, and so we use the sample mean as it estimate:
$$\widehat{\theta}_{m,k}=\frac{\sum_{i=1}^m X_{m,k}\xi_{j,k}}{N_{m,k}},\;
k=1,\ldots, K. $$
In practices, if necessary, we can add $\alpha> 0$ in the numerator and $\beta> 0$ in the denominator to avoid
the nonsense case of $0/0$, or to use prior information to estimate the parameters. Assume $\ep |\xi_{m,k}|^{2+\delta}<\infty$
and write $\sigma_k^2=\Var\{\xi_{m,k}\}$. Write $Q_{m,k}=\sum_{i=1}^mX_{m,k}(\xi_{m,k}-\theta_k)$, then $Q_{m,k}$ is a martingale
and $\widehat{\theta}_{m,k}-\theta_k=Q_{m,k}/N_{m,k}$. By the LIL (Theorem B), we have
\begin{equation}\label{eqLILofQ}
Q_{m,k}=O(\sqrt{n\log\log n}) \;\; a.s.
\end{equation}
We first show the consistency of $\bm N_n/n$.
If let $\Cal A_m=\sigma(\bm X_1,\ldots, \bm X_{m+1}$, $\bm \xi_1$, $\ldots,\bm \xi_m)$, then
$ \sum_{i=1}^m \ep[ (X_{i,k}(\xi_{i,k}-\theta_k))^2|\Cal A_{i-1}]=\sigma_k^2 N_{m,k}$.
By Theorem A (b), it follows that
\begin{equation}\label{eqthetaconver}\widehat{\theta}_{m,k}\to \theta_k\;\; a.s. \;\;\text{on the event } \; \{N_{m,k}\to \infty\}.
\end{equation}
On the event $\{N_{m,k}<\infty\}$, $\widehat{\theta}_{m,k}$ will fix
to a value eventually. In either case, $\widehat{\theta}_{m,k}$ has
a limit $\widetilde{\theta}_k$ in the parameter space,
$k=1,\cdots,K$. By the continuity of $\bm\rho(\cdot)$,
$\widehat{\bm\rho}_{m-1} \to \bm\rho(\widetilde{\theta}_1,\ldots,
\widetilde{\theta}_K ) :=\widetilde{\bm v}\in (0,1)^{\otimes K}$
a.s. Notice that the minimum of $\sum_j y_j(y_j/x_j)^{\gamma}$ over
$\sum_jx_j=1$ and $x_j\ge 0$ is $\sum_jy_j$. It is easily seen that
$g_k(\bm x,\bm y)\le y_k(y_k/x_k)^{\gamma}< y_k$ if $x_k> y_k$ and
$\sum_jx_j=\sum_jy_j=1$. So, $P_{m,k}\le \widehat{\rho}_{m-1,k}$ if
$N_{m-1,k}/(m-1)>\widehat{\rho}_{m-1,k}$. Denote
$M_{n,k}=\sum_{m=1}^n(X_{m,k}-\ep[X_{m,k}|\Cal F_{m-1}])$. Let
$S_n=\max\{m\ge MK+1: N_{m-1,k}\le (m-1)\widehat{\rho}_{m-1,k}\}$
and $\max\{\emptyset\}=MK$. Then
\begin{align}\label{eqconsistency1}
&N_{n,k}= N_{S_n,k}+\sum_{m=S_n+1}^n (X_{m,k}-\ep[X_{m,k}|\Cal F_{m-1}])+\sum_{m=S_n+1}^n P_{m,k}
\nonumber\\
\le & 1+N_{S_n-1,k}+M_{n,k}-M_{S_n,k}+\sum_{m=S_n}^{n-1}\widehat{\rho}_{m,k}
\nonumber\\
\le & 1+ N_{MK-1,k}+(M_{n,k}-M_{S_n,k})+(S_n-1)\widehat{\rho}_{S_n-1,k}+\sum_{m=S_n}^{n-1}\widehat{\rho}_{m,k}.
\end{align}
Notice that $|M_{n,k}-M_{S_n,k}|\le \max_{m\le
n}|M_{m,k}|=O(\sqrt{n\log\log n})$ a.s. by the LIL, and
$\widehat{\rho}_{m,k}\to \widetilde{v}_k$. We conclude that
$\limsup_{n\to \infty} N_{n,k}/n\le \widetilde{v}_k$ a.s.,
$k=1,\ldots, K$. From the fact that $\sum_{k=1}^n
N_{n,k}/n=\sum_{k=1}^n \widetilde{v}_k=1$, we conclude that $
N_{n,k}/n\to \widetilde{v}_k\in (0,1)$ a.s., which implies that
$N_{n,k}\to \infty$,
 $k=1,\ldots, K$. Hence $\widetilde{\theta}_k$ and $\theta_k$ ($\widetilde{\bm v}$ and $\bm v$)
 must be identical by (\ref{eqthetaconver}).
We have proved the consistency of $\bm N_n/n$. Further, according to
(\ref{eqLILofQ}), we have
$\widehat{\bm \theta}_m-\bm \theta=O\big(\sqrt{\log\log m}/\sqrt{m}\big)$ a.s.,
and then $\widehat{\bm \rho}_m-\bm v=O\big(\sqrt{\log\log m}/\sqrt{m}\big)$ a.s., which together with
 (\ref{eqconsistency1}), yields
 \begin{align*}
N_{n,k}&-nv_k
 \le  O(\sqrt{n\log\log n})+O(\sqrt{S_n\log\log S_n})\\
 &+
 \sum_{m=S_n}^{n-1}O\big(\sqrt{\log\log m}/\sqrt{m}\big)
 = O(\sqrt{n\log\log n})\;\; a.s.
 \end{align*}
By  the fact that $\sum_{k=1}^n N_{n,k}/n=\sum_{k=1}^n
\widetilde{v}_k=1$ again, we conclude that
\begin{equation}
\label{eqLILofN} \bm N_n-n \bm
v=O(\sqrt{n\log\log n})\;\; a.s.
\end{equation}

Now, we begin the proof of the asymptotic normality. Write $\bm
M_n=(M_{n,1},\ldots, M_{n,K})$ and $\bm Q_n=(Q_{n,1},\ldots,
Q_{n,K})$. Then by (\ref{eqLILofQ}) and (\ref{eqLILofN}),
$$ \widehat{\theta}_{m,k}-\theta_k=\frac{Q_{m,k}}{m v_k}+O\bigg(\frac{\log\log m}{m}\bigg)\;\; a.s.  $$
It is easily seen that $\partial \bm g/\partial \bm x|_{\bm x=\bm y}=-\gamma (\bm I-\bm 1^{\prime}\bm v)$
and $\partial \bm g/\partial \bm y|_{\bm x=\bm y}=(\gamma+1) (\bm I-\bm 1^{\prime}\bm v)$. By the Taylor
formula, we have
\begin{align*}
& \bm P_m-\bm v=\bm \rho\big(\frac{\bm N_{m-1}}{m-1}, \widehat{\bm \rho}_{m-1}\big)-\bm \rho(\bm v,\bm v)\\
=&-\gamma\big(\frac{\bm N_{m-1}}{m-1}-\bm v\big)(\bm I-\bm 1^{\prime}\bm v)
+(\gamma+1)\big(\widehat{\bm \theta}_{m-1}-\bm
\theta\big)\frac{\partial \bm \rho}{\partial \bm \theta}(\bm I-\bm 1^{\prime}\bm v)\\
&+O\big(\big\|\frac{\bm N_{m-1}}{m-1}-\bm v\big\|^2\big)
+\big(\big\|\widehat{\bm \theta}_{m-1}-\bm
\theta\big\|^2\big)\\
=&-\gamma\big(\frac{\bm N_{m-1}}{m-1}-\bm v\big)
+(\gamma+1)\big(\widehat{\bm \theta}_{m-1}-\bm
\theta\big)\frac{\partial \bm \rho}{\partial \bm \theta}
+O\bigg( \frac{\log\log m}{m}\bigg)\\
=&-\gamma\big(\frac{\bm N_{m-1}}{m-1}-\bm v\big)
+(\gamma+1)\frac{\bm Q_{m-1}}{(m-1)}diag(\frac{1}{\bm v})\frac{\partial \bm \rho}{\partial \bm \theta}
+O\bigg( \frac{\log\log m}{m}\bigg)\\
=&O\bigg( \sqrt{\frac{\log\log m}{m}}\bigg).
\end{align*}
Hence
\begin{align*}
&\bm N_n-n\bm v=\sum_{m=1}^n (\bm X_m-\ep[\bm X_m|\Cal F_{m-1}])+\sum_{m=1}^n (\bm P_m-\bm v)\\
=& \bm M_n-\gamma \sum_{m=1}^{n-1}\frac{\bm N_m-m\bm v}{m}+
(\gamma+1)\sum_{m=1}^{n-1}\frac{\bm Q_m}{m}diag(\frac{1}{\bm v})\frac{\partial \bm \rho}{\partial \bm \theta}
+o((\log n)^2) \; \; a.s.
\end{align*}
On the other hand, it is easily checked that,  for the martingale $(\bm M_n,\bm Q_n)$ we have
\begin{align*}
&\sum_{m=1}^n\ep[(\Delta \bm M_m)^{\prime}(\Delta \bm M_m)|\Cal F_{m-1}]
=\sum_{m=1}^n(diag(\bm P_m)-\bm P_m^{\prime}\bm P_m)\\
& \quad = n(diag(\bm v)-\bm v^{\prime}\bm v)+O(\sqrt{n\log\log n})\;\; a.s.,\\
&\sum_{m=1}^n\ep[(\Delta Q_{m,k})^2|\Cal F_{m-1}]
=\sum_{m=1}^n P_{m,k}\sigma_k^2
 =n v_k \sigma_k^2+O(\sqrt{n\log\log n})\; a.s.,
\end{align*}
$\ep[(\Delta \bm M_m)^{\prime}(\Delta \bm Q_m)|\Cal F_{m-1}]=\bm 0$ and
$\ep[(\Delta  Q_{m,k}Q_{m,j}|\Cal F_{m-1}]=0$, $k\ne j$.
So, by the strong approximation (Theorem E), there are two independent
 $K$-dimensional standard Brownian motions
$\bm W(t)$ and $\bm B(t)$ such that for some $\kappa>0$,
$ \bm M_n=\bm W(n)\bm \Sigma_1^{1/2}+o(n^{1/2-\kappa})$ a.s.  and
$Q_{n,k}=B_k(n)\sigma_k\sqrt{v_k}+o(n^{1/2-\kappa})$  a.s., $ k=1,\ldots, K,$
where $\bm \Sigma_1=diag(\bm v)-\bm v^{\prime}\bm v$. It follows that
\begin{align*}\bm N_n-n\bm v=&
\bm W(n)\bm \Sigma_1^{1/2}-\gamma \sum_{m=1}^{n-1}\frac{\bm N_m-m\bm v}{m}+
(\gamma+1)\sum_{m=1}^{n-1}\frac{\bm B(m)}{m}\bm\Sigma_{\bm\rho}^{1/2}\\
&+o(n^{1/2-\kappa}) \; \; a.s.
\end{align*}
Hence,
$\bm N_n-n\bm v=\bm G(n)+ +o(n^{1/2-\kappa}) \; \; a.s.,$
where
$$ \bm G(t)=t^{-\gamma}\int_0^t x^{\gamma} d \bm W(x) \bm \Sigma_1^{1/2}+
(\gamma+1)t^{-\gamma}\int_0^t x^{\gamma-1} \bm B(x) dx \bm \Sigma_{\bm\rho}^{1/2}$$
is the solution of the equation
$$ \bm G(t)=\bm W(t)\bm \Sigma_1^{1/2}-\gamma\int_0^t \frac{\bm G(x)}{x}dx +
(\gamma+1)\int_0^t \frac{\bm B(x)}{x}dx \bm \Sigma_{\bm\rho}^{1/2}. $$
It follows that $\sqrt{n}(\bm N_n/n-\bm v)\overset{\mathscr D}\to N(\bm 0,\bm\Sigma)$ with
\begin{align*}
 &\bm\Sigma=\Var\{\bm G(1)\}\\
 =&\bm\Sigma_1\int_0^1x^{2\gamma}dx+(\gamma+1)^2
 \bm\Sigma_{\bm\rho}\int_0^1\int_0^1x^{\gamma-1}y^{\gamma-1}
(x\wedge y) dxdy\\
 =& \frac{1}{2\gamma+1}\bm\Sigma_1+\frac{2(\gamma+1)}{2\gamma+1}\bm\Sigma_{\bm\rho}.
 \end{align*}
The proof is now completed. $\Box$


\section{Delayed responses}\label{section7}
\setcounter{equation}{0}

In practices, the outcomes in clinical trials are not available
immediately prior to the treatment allocation of the next subject.
The estimating of the  parameters, and the updating of the urn when
using urn models, can only be processed according to observed
responses. The effect of the delay of treatment results is fist
studied in theory by  Bai, Hu and Rosenberger (2002) for the urn
compositions in an urn model design with discrete responses. After
that, Hu and Zhang (2004b), Zhang, Chan, Chueng and Hu (2007), Sun,
Cheung and Zhang (2007) and Zhang, Hu, Chueng and Chan (2006a) have
shown that the delay machine does not effect the asymptotic
properties of the sample allocation proportions for many adaptive
designs if the delay degree decays with a power rate. The basic
reason is that the total delayed responses
 is a high order of square root of the sample size when the delay degree decays with a power rate.

To describe the delay machine, we let $t_m$ be the entry time of the
$m$-th subject, where $t_m$ is an increasing sequence of random
variables. Assume that $\{t_{m+1}-t_m\}$ is a sequence of
independent random variables.  The response time of the $m$-th
subject on treatment $k$ is denoted by $r_m(k)$. Suppose $\{r_m(k);
m\ge 1\}$ are sequence of independent random variables, $k=1,\ldots,
K$. Further, assume that $\{t_{m+1}-t_m, r_m(k);k=1,\cdots,K, m\ge n\}$
is independent of the assignments $\bm X_1,\ldots,\bm X_n$.

\begin{assumption}\label{Asp1}
Let $\delta_k(m,l)=I\{r_m(k)\ge t_{m+l}-t_m\}$ be an indicator
function that takes the value $1$ if the outcome of the $m$-th
subject on treatment $k$ occurs after at least another $l$ subjects
arrive, and $0$ otherwise. Suppose for some constants $C>0$ and
 $\gamma\ge 2$,
$$ \mu_k(m,l)=\pr\{\delta_k(m,n)=1\}\le C l^{-\gamma},\;\;  m,l=1,2,\ldots, k=1,\ldots, K.
$$
\end{assumption}
This assumption is easily satisfied. A practical approach is to
assume that the entry mechanism generates a Poisson process and the
delay time has an exponential distribution in which both
$\{r_m(k)\}$ and   $\{t_{m+1}-t_m\}$  are sequences of  i.i.d.
exponential random variables with means $\lambda_k>0$ and
$\lambda_0>0$, respectively. This approach is common in clinical
studies and the probability $\mu_k(m,l)$ is
$\big(\lambda_k/(\lambda_0+\lambda_k)\big)^l$.

Let $S_{m,k}^{obs}$ (resp. $N_{m,k}^{obs}$) be the sum (resp. the
number) of the outcomes on treatment $k$ observed prior to the
$(m+1)$-th assignment, and  $S_{m,k}$ (resp. $N_{m,k}$) be the sum
(resp. the number) of all the outcomes of those being assigned to
treatment $k$ in the first $m$ subjects, $k=1,\ldots,K.$

\begin{theorem}\label{thedelay} Suppose Assumptions
\ref{Asp1} is  satisfied, and the responses on each treatment are
i.i.d. random variables having finite $(2+\delta)$-th moments. Then
for some $0<\delta_0<\frac{1}{2}-\frac{1}{2+\delta}$, we have
\begin{equation}\label{eqlem1.3}
R_{n,k}=:S_{m,k}-S_{m,k}^{obs}=o(n^{1/2-\delta_0}) \quad a.s.
\end{equation}
and $N_{n,k}-N_{n,k}^{obs}=o(n^{1/2-\delta_0})$ a.s., $k=1,\ldots, K$.
\end{theorem}

{\bf Proof.} Let $I_k(m,l)$ be the indicator function, which takes
value $1$ if the outcome $\xi_{m,k}$  on treatment $k$ of the $m$-th
subject
 occurs after the $(m+l)$-th
assignment and before the $(m+l+1)$-th assignment, $k=1,\ldots, K$.
For given $m$ and $l$, if $I_k(m,l)=1$, we observe a response
$X_{m,k}\xi_{m,k}$.  Hence, between the $m$-th assignment and the
$(m+1)$-th assignment, the observed outcomes are
$I_k(m,0)X_{m,k}\xi_{m,k}$, $I_k(m-1,1)X_{m-1,k}\xi_{m-1,k}$,
$\ldots$, $I_k(1,m-1)X_{1,k}\xi_{1,k}$, $k=1,\ldots, K$, and the sum
of those on treatment $k$ is $
\sum_{l=0}^{m-1}I_k(m-l,l)X_{m-l,k}\xi_{m-l,k}=\sum_{l=1}^m
I_k(l,m-l)X_{l,k}\xi_{l,k}$. Hence,
$$S_{m,k}^{obs}=\sum_{m=1}^n\sum_{l=1}^m
I_k(l,m-l)X_{l,k}\xi_{l,k}=\sum_{m=1}^n \sum_{j=m}^nX_{m,k}
\xi_{m,k}I_k(m, j-m).
$$
On the other hand, it is obvious that
$S_{m,k}=\sum_{m=1}^nX_{m,k}\xi_{m,k}$ $=\sum_{m=1}^n $ $\sum_{j=m}^{\infty}X_{m,k}\xi_{m,k}I_k(m,
j-m)$. It  follows that
$$ R_{n,k}=\sum_{m=1}^n\sum_{j=n+1}^{\infty}X_{m,k}\xi_{m,k}I_k(m,
j-m). $$
 Let $0<\phi<1/2$ be a number whose value will be
specified later. Write $l_n=[n^{\phi}]$. Then
\begin{align}\label{eqproofofRnk}
&|R_{n,k}|=\Big|\Big(\sum_{m=n-l_n+1}^n+\sum_{m=1}^{n-l_n}
\Big)
\sum_{j=n-m+1}^{\infty}I_k(m,j)X_{m,k}\xi_{m,k}\Big|\nonumber\\
\le & \sum_{m=n-l_n+1}^n X_{m,k}|\xi_{m,k}|+\sum_{m=1}^n
\delta_k(m,l_m)X_{m,k}|\xi_{m,k}|\nonumber\\
\le &Cn^{\phi}+\sum_{m=n-l_n+1}^n
 \big(|\xi_{m,k}|-\ep[|\xi_{m,k}|]\big)+\sum_{m=1}^n X_{m,k}\ep[\delta_k(m,l_m)|\xi_{m,k}|]
\nonumber\\
&+\sum_{m=1}^n
X_{m,k}\big(\delta_k(m,l_m)|\xi_{m,k}|-\ep[\delta_k(m,l_m)|\xi_{m,k}|]\big).
\end{align}
For  $\sum_{m=1}^n \big(|\xi_{m,k}|-\ep[|\xi_{m,k}|]\big)$, due to
Theorems 1.2.1 and 2.6.6 of Cs\"org\H o and R\'ev\'esz (1981),
\begin{align*}
\sum_{m=n-l_n+1}^n \big(|\xi_{m,k}|-\ep[|\xi_{m,k}|]\big)=&
O(\sqrt{l_n\log n})+o(n^{\frac{1}{2+\delta}})\\
=& o(n^{\phi})+o(n^{\frac{1}{2+\delta}}) \;\; a.s.
\end{align*}
 For the martingale $\sum_{m=1}^n
X_{m,k}\big(\delta_k(m,l_m)|\xi_{m,k}|-\ep[\delta_k(m,l_m)|\xi_{m,k}|]\big)$,
notice \begin{align*} &\ep\Big(
X_{m,k}\big(\delta_k(m,l_m)|\xi_{m,k}|-\ep[\delta_k(m,l_m)|\xi_{m,k}|]\big)\Big)^2\\
\le & \ep\big( \delta_k(m,l_m)|\xi_{m,k}|\big)^2 \le
 \big(\ep
\delta_k(m,l_m)\big)^{\frac{\delta}{2+\delta}}\big(\ep|\xi_{m,k}|^{2+\delta}\big)^{\frac{2}{2+\delta}}
\le  m^{-\frac{\phi\gamma\delta}{2+\delta}}.
\end{align*}
due to Assumptions \ref{Asp1} and the assumption of finite $(2+\delta)$-th moments. So, by the LLN (Theorem A (a)),
$\sum_{m=1}^n
X_{m,k}\big(\delta_k(m,l_m)|\xi_{m,k}|-\ep[\delta_k(m,l_m)|\xi_{m,k}|]\big)
=o(n^{\frac{1}{2}-\frac{\phi\gamma\delta}{2(2+\delta)}}\log n)$
a.s.
 Also
\begin{align*} &\sum_{m=1}^nX_{m,k}\ep[\delta_k(m,l_m)|D_{m,k}|] \le
\sum_{m=1}^n\big(\ep
\delta_k(m,l_m)\big)^{\frac{1+\delta}{2+\delta}}\big(\ep|D_{m,k}|^{2+\delta}\big)^{\frac{1}{2+\delta}}\\
& \quad \le  C\sum_{m=1}^n m^{-\frac{\phi\gamma(1+\delta)}{2+\delta}}\le
C n^{1-\frac{\phi\gamma(1+\delta)}{2+\delta}}.
\end{align*}
Combining the above arguments yields
$$R_{n,k}=O(n^{\phi})+O(n^{1-\frac{\phi\gamma(1+\delta)}{2+\delta}})+o(n^{\frac{1}{2+\delta}})+
o(n^{\frac{1}{2}-\frac{\phi\gamma\delta}{2(2+\delta)}}\log n)\;\;
a.s.
$$
Choosing $\phi=\frac{2+\delta}{1+(1+\delta)\gamma}$ and
$0<\delta_0<\min\{\frac{1}{2}-\phi,\frac{1}{2}-\frac{1}{2+\delta},
\frac{\phi\gamma\delta}{2(2+\delta)}\}$ yields (\ref{eqlem1.3}).
The proof of $N_{n,k}-N_{n,k}^{obs}=o(n^{1/2-\delta_0})$ a.s. is similar. $\Box$

\section{Variability, power and asymptotic best adaptive designs}\label{section6}
\setcounter{equation}{0}

 The variability  of the sample proportion
$N_{n,k}/n$ is an important quantity, which   measures the distance
between $N_{n,k}/n$ and its limit $v_k$. The smaller is the
variability, the smaller is the probability that there is large bias
between $ N_{n,k}/n$ and $v_k$. When using a adaptive design with
high variability, a clinical might result in assigning more subjects
to the inferior treatment making the allocation even less ethical
than equal allocation. Also, a trial with high variability might
result in assigning only a few subjects to one of the the treatments
decreasing the efficiency in the test or the estimation  of
parameters. According to (\ref{eqvarofRPW}),  the asymptotic
variability of the RPW rule is very high unless both treatments have
low success rates. It is extremely high when $q_1+q_2$ is close to
$1/2$, and it is showed that $n\Var\{N_{n,1}/n\}\to \infty$ when
$q_1+q_2<1/2$. The RPW rule used in ECMO, 1985, trial assigned only
one patient to the less successful control therapy (see Royall,
1991, for discussion). The  relationship among the power, the target
allocation and the variability of the designs is fist revealed by Hu
and Rosenberger (2003) in theory, though simulation studies  had
indicated there is strong relationship among these quantities. Hu
and Rosenberger (2003) proved that the average power of a
statistical test of the difference of distribution parameters is a
decreasing function of the variability of the designs. In Section
\ref{section5}, we have found that the asymptotic variability of
DBCD is a decreasing function of the parameter $\gamma$. When
$\gamma\to \infty$, the variability tends to its  minimums
$$\bm \Sigma_{\rm \rho}=\left(\frac{\partial \bm\rho}{\partial \bm\theta}\right)^{\prime}
diag\left((v_1 I_1(\theta_1))^{-1},\ldots,v_K I_1(\theta_K))^{-1}\right)
\frac{\partial \bm\rho}{\partial \bm\theta}.
$$
Hu,  Rosenberger and Zhang (2006) proved that this limit is the
lower bound of the asymptotic variability of adaptive designs among
all adaptive designs which have the same  limiting proportion.

Assume the following regularity conditions:
\begin{enumerate}
  \item The parameter space $\bm\Theta_k$ is an open subset in $\mathscr R^d$, $j=1,\ldots,K$;
  \item The distributions of outcomes $f_1(\cdot|\bm \theta_1)$, $\ldots$, $f_K(\cdot|\bm \theta_K)$
  follow an exponential family;
  \item For the limiting allocation proportion $\bm\rho(\bm \theta)=(\rho_1(\bm \theta),\ldots,\rho_K(\bm \theta)
  \in (0,1)^{\otimes K}$,
  $$ \frac{\bm N_{n,j}}{n}\to \rho_j(\bm \theta) \;\; a.s. \;\; j=1,\ldots, K;$$
  \item For a positive definite matrix $\bm V(\bm \theta)$,
  $$ \sqrt{n}\left(\frac{\bm N_n}{n}-\bm \rho(\bm \theta)\right)\overset{\mathscr D}\to N(\bm 0, \bm V(\bm\theta). $$
\end{enumerate}
\begin{theorem} Under regularity conditions 1-4, there exists a $\bm \Theta_0\subset \bm \Theta=\bm\Theta_1\otimes\cdots
\otimes\bm \Theta_K$ with Lebesgue measure $0$ such that for every $\bm \theta\in \bm\Theta-\bm\Theta_0$,
$$ \bm V(\bm\theta)\ge \left(\frac{\partial \rho(\bm\theta)}{\partial \bm\theta}\right)^{\prime}
\bm I^{-1}(\bm \theta)\frac{\partial \rho(\bm\theta)}{\partial \bm\theta}:=\bm\Sigma_{LB}, $$
where $\bm I(\bm\theta)=diag\big(\rho_1(\bm\theta)\bm I_1(\bm\theta_1),\ldots,\rho_K(\bm\theta) \bm I_K(\bm\theta_K)\big)$
and $\bm I_k(\bm\theta_k)$ is the Fisher information for a single observation on treatment $k=1,\ldots, K$.
\end{theorem}

 We refer to a adaptive design that attains the lower bound as
 {\em asymptotic best }
 for that particular allocation $\bm \rho(\bm\theta)$.
 Table \ref{table1} gives the lower bounds of
the asymptotic variabilitites ($\sigma_{LB}^2$) for urn proportion
(UP), Roenberger, et al's proportion (RP) and Neyman proportion (NP)
in a binary response clinical trial for two-treatments.
\begin{table}[h!]
  \centering
  \begin{tabular}{|c|c|c|}
    \hline
       & $\rho$                     &  $\sigma_{LB}^2$ \\
  \hline
    UP & $\frac{q_2}{q_1+q_2}$      & $\frac{q_1q_2(p_1+p_2)}{(q_1+q_2)^3}$
    \\[5pt]
    RP & $\frac{\sqrt{p_1}}{\sqrt{p_1}+\sqrt{p_2}}$ & $\frac{1}{4 (\sqrt{p_1}+\sqrt{p_2})^3}\left(\frac{p_2q_1}{\sqrt{p_1}}+\frac{p_1q_2}{\sqrt{p_2}}\right)$
    \\ [5pt]
    NP & $\frac{\sqrt{p_1q_1}}{\sqrt{p_1q_1}+\sqrt{p_2q_2}}$ & $\frac{1}{4(\sqrt{p_1q_1}+\sqrt{p_2q_2})^3}
\left(\frac{p_2q_2(1-2p_1)^2}{\sqrt{p_1q_1}}
+\frac{p_1q_1(1-2p_2)^2}{\sqrt{p_2q_2}}\right)$ \\[5pt]
    \hline
  \end{tabular}
  \caption{Lower bounds of the asymptotic variabilities for urn proportion,
  Rosenberger et al's proportion and Neyman proportion.}
  \label{table1}
\end{table}
It is interesting to note that the Zelen (1969)'s deterministic
design (PW rule) is
 asymptotic best among all procedures with limiting allocation proportion $q_2/(q_1+q_2)$, the urn proportion.
The RPW rule has the same limiting  proportion, but it is not
asymptotic best because $\sigma^2_{RPW}>\sigma_{LB}^2$. The DL rule
proposed by Ivanova (2003) is a random procedure having the same
urn proportion   and the asymptotic variability
$\sigma_{DL}^2=q_1q_2(p_1+p_2)/(q_1+q_2)^3$. So, the DL rule  is
asymptotic best. However, both the RPW rule and the DL rule can only
target this particular proportion, which is not
optimal in any formal sense, and can only be used for binary
responses (c.f., Hu and Rosenberger, 2003).

DBCD is also not asymptotic best (except $\gamma=\infty$). But, it
can target any desired allocation and can be used for general
responses, for example, continuous response. Zhang, Hu and Cheung
(2006) proposed an urn model, the SEU model, which can target any
pre-specified allocation proportion and can be used for general
responses. The drop-the-loser rule has also been generalized to GDL
model, a kind of urn model with immigration,  by Zhang, Chan, Cheung
and Hu (2007) and Sun, Chueng and Zhang (2007), by using the
estimators of unknown distribution parameters, such that it can
target any pre-specified allocation proportion and can be used for
general responses. For a general pre-specified  allocation
proportion $\bm\rho=\bm\rho(\bm\theta)$ in a $K$-treatment trial,
the asymptotic variance-covariance matrices of SEU, GDL and DBCD are
given in Table \ref{table2}.
\begin{table}[h!]
  \centering\small
  \begin{tabular}{|c|c|c|c|}
    \hline
    Model & SEU & GDL  & DBCD \\
    \hline
    Variability  & $diag(\bm\rho)-\bm\rho^{\prime}\bm\rho+6\bm\Sigma_{LB}$ &$2\bm\Sigma_{LB}$   &  $\frac{diag(\bm\rho)-\bm\rho^{\prime}\bm\rho}{1+2\gamma}+
    \frac{2+2\gamma}{1+2\gamma}\bm\Sigma_{LB}$
    \\[5pt]
    \hline
  \end{tabular}
  \caption{The asymptotic variability  of SEU, GDL and DBCD for a same limiting allocation }\label{table2}
\end{table}
Among these models, the DBCD can approach the lower bound for large
values of $\gamma$. However, the procedure becomes more
deterministic as $\gamma$ becomes larger, and hence careful tuning
of $\gamma$ must be done to counter the trade-off between the
randomness and variability. The use of DBCD with $\gamma=2$ was
strongly recommended in Hu and Rosenberer (2003) for binary response
trials with two treatments, according to the simulation study.

Very recently, we have found a fully randomized biased coin design
(RBCD) for two-treatment clinical trails,
 a kind of DBCD, which preserves randomization,
attains the lower bound, and can target any allocation. Whether or
not an urn model can be defined to have these properties is still an
open problem.
In the RBCD, instead of using a continuous allocation
function $g(x,y)$, we use a discrete function: $g(x,y)=\alpha y$ if
$x>y$, $y$ if $x=y$ and $1-\alpha(1-y)$ if $x<y$, where
$0<\alpha<1$.
 That is, the $m$-th subject is allocated to treatment $1$ with a
 probability
\begin{align*}&\pr(X_{m,1}=1|\Cal F_{m-1}) =P_{m,1}\\
=&\begin{cases}\alpha \rho(\widehat{\bm\theta}_{m-1}), & \text{ if }
N_{m-1,1}/(m-1)>\rho(\widehat{\bm\theta}_{m-1}),\\
\rho(\widehat{\bm\theta}_{m-1}), & \text{ if }
N_{m-1,1}/(m-1)=\rho(\widehat{\bm\theta}_{m-1}),\\
1-\alpha(1- \rho(\widehat{\bm\theta}_{m-1})), & \text{ if }
N_{m-1,1}/(m-1)<\rho(\widehat{\bm\theta}_{m-1}).
\end{cases}
\end{align*}
When $\rho(\bm\theta)\equiv 1/2$, the RBCD is just the Efron
(1971)'s biased coin design. The following theorem gives  the
asymptotic results for the RBCD, the proof of which will not be
presented here.
\begin{theorem}
 Suppose $\rho(\cdot,\cdot)$ is a continuous function on the parameter space and
 twice differentiable at the true value of the parameter
$\bm\theta=(\theta_1,\theta_2)$,
 the distributions $f_1(\cdot|\theta_1)$ and $f_2(\cdot|\theta_2)$
 of the responses follow an exponential family. Then
$N_{n,1}/n-\rho(\bm\theta)=O(\sqrt{\log\log n/n})$  a.s. and
$ \sqrt{n}\left(N_{n,1}/n-\rho(\bm\theta)\right)\overset{D}\to
N(0,\sigma_{\rho}^2),
$
where
$$\sigma^2_{\rho}= \left(\partial \rho(\bm\theta)/\partial
\theta_1\right)^2(\rho(\bm\theta)I_1(\theta_1))^{-1}+ \left(\partial
\rho(\bm\theta)/\partial
\theta_2\right)^2((1-\rho(\bm\theta))I_2(\theta_2))^{-1}. $$
Further, $ \max_{m\le n}|N_{m1}-m\rho(\bm\theta)-\sigma_{\rho}
W(m)|=o(\sqrt{n})$  in probability, where $W(t)$ is a standard
Brownian motion.
\end{theorem}

\section{Discussion}\label{section8}

In this paper, we have discussed several classes of adaptive
designs. The play-the-winner  rule  is the simplest procedure and
has small variability. But is is too deterministic to be used in
clinical trials. The randomized play-the-winner rule and urn models
are random procedures. But their variabilities are very high.
Theoretical results, simulation studies and a real example in ECMO
trial all indicate  that the statistical test in using a adaptive
design with high variability  is not powerful. However, besides in
adaptive designs, urn models have wide applications in many areas
including biological science, random algorithm and sampling,
information science, etc. And the urn models have strong
relationship with multi-type branching processes. The study of urn
models has been of interest in a long history.   The DL rule is
 randomized procedure and has the smallest variability among all the
adaptive designs with limiting allocation proportion
$q_2/(q_1+q_2)$. But it
 can only target this particular proportion  and can only be used for binary
 responses. When it is generalized to be able to target any
 pre-specified allocation, the variability is no longer the smallest
 (c.f. Table \ref{table2}). Among the adaptive designs mentioned
 in this paper, the DBCD and RBCD are the only procedures that preserves randomization,
attains or can approach the lower bound of the variability, can
target any allocation and can be used for general discrete or
continuous responses.

The examples considered in this paper are binary response clinical
trials. In practices, the responses of clinical trials appear in
various types. For more examples and discussion, we refer to a new
book of Hu and Rosneberger (2006). In many clinical trials,
covariate information is available that has a strong influence on
the responses of patients. For instance, the efficacy of a
hypertensive drug is related to a patient's initial blood pressure
and cholesterol level, whereas the effectiveness of a cancer
treatment may depend on whether the patient is a smoker or a
non-smoker. The theory of an adaptive design in using covariate
information is much more complicated than those without covarites. A
limit success in deriving the asymptotic properties of
covariate-adjusted  adaptive designs has been achieved  by Zhang, Hu
Cheung and Chan (2007). The power study and evaluation of the
covariate-adjusted  adaptive designs are our future studies. Also,
in many clinical trials, the observed responses are usual survival
data. Though it has been shown that the delay of treatment results
does not effect the asymptotic properties in many adaptive designs,
 is is assume that the delayed responses are finally observed if
the time is long enough. It is an interesting topic of studying the
properties of the adaptive designs with missed or censored data.

Finally, in this paper, we only present the asymptotic results. It
is important to check the accuracy of the asymptotic approximations
when using the theoretical results to evaluate or compare designs.
Simulations have indicated  that in most cases these designs closely
approximate asymptotic results for a moderate sample size of
$n=100$.




\begin{thebibliography}{99} 

\bibitem{AK} {\sc Athreya, K. B. \& Karlin, S.}  (1968).
   Embedding of urn
schemes into continuous time branching processes and related limit
theorems.   {\em Ann. Math. Statist.}  {\bf 39}   1801-1817.

\bibitem{BH99} {\sc Bai, Z. D. \& Hu, F.}  (1999). Asymptotic theorem
for urn models with nonhomogeneous generating matrices. {\em
Stochastic Process. Appl.} {\bf 80}  87-101.

\bibitem{BH05} {\sc Bai, Z. D. \& Hu, F.} (2005).
 Strong consistency and asymptotic
normality for urn models.
 {\em Ann. Appl. Probab.} {\bf 15} 914--940.

\bibitem{BHR02} {\sc Bai, Z.D., Hu, F.\& Rosenberger, W.F.} (2002). Asymptotic
properties of adaptive designs for clinical trials with delayed
response. {\em Ann. Statist.} {\bf 30}, 122-139.

\bibitem{BHS} {\sc Bai, Z. D., Hu, F. \& Shen, L.} (2002).
An adaptive design for multi-arm clinical trials.
 {\em J. Multi.
Anal.} {\bf 81} 1--18.


\bibitem{BHZ} {\sc Bai, Z. D., Hu, F. \& Zhang, L.-X.}  (2002). The Gaussian approximation theorems for urn
 models and their applications. {\em Ann. Appl.  Probab.} {\bf
12} 1149--1173.

\bibitem{Connor} {\sc Connor, E. M., Sperling, R. S., Gelber, R., Kiselev,
P., Scott, G., O'Sullivan, M. J., VanDyke, R., Bey, M., Shearer, W.,
Jacobson, R. L.,  Jiminez, E., O'Neil, E., Bazin, B., Delfraissy,
J., Culnane, M., Coombs, R., Elkins, M., Moye, J., Stratton, P. \&
Balsley, J for the Pediatric AIDS Clinical Trials Group Protocol 076
Study Group} (1994). Reduction of maternal-infant transmission of
human immunodeficiency virus type 1 with zidovudine treatment. {\em
New England Journal of Medicine} {\bf 331} 1173-1180.

\bibitem{Efron} {\sc Efron, B.} (1971). Forcing a sequential
experiment to be balanced.
 {\em Biometrika} {\bf 62}  347-352.

\bibitem{Ei94} {\sc Eisele, J.} (1994).
The doubly adaptive biased coin design for sequential clinical
trials. {\em J. Statist. Plann. Inf.}, {\bf 38}: 249-262.

\bibitem{EW95} {\sc Eisele, J. \& Woodroofe, M.} (1995).
Central limit theorems for doubly adaptive biased coin designs. {\em
Ann. Statist.} {\bf 23} 234-254.


\bibitem{HallHeyde} {\sc Hall, P. \&  Heyde, C. C.} (1980). {\em
Martingale Limit Theory and its Applications}. Academic Press,
London.

\bibitem{HR03} {\sc Hu, F. \&  Rosenberger, W. F.}  (2003).
Optimality, variability, power   evaluating response-adaptive
randomization procedures for treatment comparisons.  {\em J. Amer.
Statist. Assoc.} {\bf 98} 671-678.

\bibitem{HR06}
{\sc Hu, F. \& Rosenberger, W. F.}  (2006).  {\em The Theory of
Response-Adaptive Randomization in Clinical Trials}. John Wiley and
Sons. Wiley Series in Probability and Statistics.



\bibitem{HRZ} {\sc Hu, F., Rosenberger, W. F. \& Zhang, L.-X.}  (2006).
 Asymptotically best
response-adaptive randomization procedures. {\em J. Statist. Plann.
Inf.}  {\bf 136} 1911--1922.


\bibitem{HZ01} {\sc Hu, F. \& Zhang, L.-X.} (2001). The weak and strong
 invariance for generalized Friedman's urn model. {\em Manuscript}.
 \underline{http://www.math.zju.edu.cn/zlx/mypapers2/URNHDIM3.pdf}

\bibitem{HZ04a} {\sc Hu, F. \& Zhang, L.-X.}  (2004a).  Asymptotic properties of doubly
adaptive biased coin designs for multi-treatment clinical trials.
{\em Ann. Statist.} {\bf 32} 268--301.

\bibitem{HZ04b} {\sc Hu, F. \& Zhang, L.-X.}  (2004b). The asymptotic normality of urn
models for clinical trials with delayed response. {\em Bernoulli}
{\bf  10 } 447-463.

\bibitem{Iv03} {\sc Ivanova, A. V.}  (2003).  A play-the-winner type urn model
with reduced variability.  {\em Metrika} {\bf 58} 1--13.

\bibitem{Iv06} {\sc Ivanova, A. V.} (2006). Urn designs with immigration: Useful connection with
continuous time stochastic processes. {\em  J. Statist. Plann. Inf.}
{\bf 136}  1836-1844.

\bibitem{IR}
{\sc Ivanova, A. V. \& Rosenberger, W. F.} (2003). A comparison of
urn designs for randomized clinical trials of $K>2$ treatments. {\em
J. Biopharmaceutial Statistics} {\bf 10} 93--107.

\bibitem{Janson}{\sc Janson, S.} (2004). Functional limit theorems for multitype
branching processes and generalized P\'olya urns. {\em Stochastic
Process. Appl.}.  {\bf 110} 177--245.

\bibitem{JT} {\sc Jennison,  C. \& Turnbull, B. W.} (2000).
{\em Group Sequential Methods with Applications to Clinical Trials}.
Chapman and Hall/CRC, Boca Raton, FL.

\bibitem{LBCH} {\sc  Lin, Z. Y., Bai, Z. D., Chen, Y. M. \& Hu, F. }
(2003). Adaptive designs based on Markov chains (I) (II). {\em J.
Biomath.} {\bf 18} 597-609 (in Chinese).


\bibitem{LZCC} {\sc Lin, Z. Y., Zhang, L. X., Chueng, S. H. \& Chan, W.
S. } (2005). Strong approximation of a Markov chain. {\em Chinese
Ann. Math.} {\bf 26A} 241-250 (in Chinese); {\em Chinese J.
Contemporary Math.} {\bf 26} (2006) 283-290 (in English).

\bibitem{Robbins} {\sc Robbins,} H. (1952). Some aspects of the
sequential design of experiments. {\em Bull. Amer. Math. Soc.} {\bf
58} 527-535.

\bibitem{Rosen} {\sc Rosenberger, W. F.} (1999). Randomized
play-the-winner clinical trials: review and recommendation. {\em
Controlled Clinical Trials}. {\bf 2}. 328-342.

\bibitem{RFD} {\sc Rosenberger, W. F.,  Flournoy, N. \&  Durham, S. D.} (1997). Asymptotic normality of maximum
likelihood estimators from multiparameter response-driven design.
{\em J. Statist. Plann. Inf.} {\bf 60} 69-76.

\bibitem{RH04} {\sc Rosenberger, W. F. \&  Hu, F.} (2004).
Maximizing power and minimizing treatment failures. {\em Clinical
Trials}  {\bf 1} 141-147.

\bibitem{RL} {\sc Rosenberger, W. F. \& Lachin, J. M.} (1995).
{\em Randomization in Clinical Trials: Theory and Practice}. Wiley,
New York.

\bibitem{RSIHR} {\sc Rosenberger, W. F., Stallard, N., Ivanova, A., Harper, C.
N., \& Ricks, M. L.}  (2001).  Optimal adaptive designs for binary
response trials.  {\em Biometrics} {\bf 57}  909--913.

\bibitem{SVA01} {\sc Rosenberger, W. F.,  Vidyashankar,  A. N.
\& Agarwal,} D. K. (2001). Covariate-adjusted response-adaptive
designs for binary response. {\em J. Biopharmaceutial Statist.} {\bf
11} 227-236.

\bibitem{Royall} {\sc Royall, R. M.} (1991). Ethics and statistics in
randomized clinical trials (Disc:p63-88). {\em Statistical Science}
{\bf 6} 52-62.

\bibitem{SR}
{\sc Smythe, R. T. \& Rosenberger, W. F.} (1995). Play-the-winner,
generalized Polya urns, and Markov branching processes.  In {\em
Adaptive Designs} ({\sc Flournoy, N. \& Rosenberger, W. F., eds.}).
Hayward:  Institute of Mathematical Statistics, 13--22.

\bibitem{Stout74} Stout, W. F. (1974). {\em Almost sure convergence}.
 Academic Press, New York.

\bibitem{SCZ} {\sc Sun, R., Cheung, S. H. \&  Zhang, L.-X.} (2007).
 A generalized drop-the-loser rule for multi-treatment clinical trials.
 {\em J. Statist. Plann.
Inf.} {\bf 137} in press.


\bibitem{Thom} {\sc Thompson, W. R.} (1933). On the likelihood that
one unknown probability exceeds another in view of the evidence of
the two samples. {\em Biometrika} {\bf 25} 275-294.

\bibitem{Wei79} {\sc Wei, L. J.}  (1979).  The generalized Polya's urn
design for sequential medical trials.  {\em Ann. Statist.}, {\bf 7}
291--296.

\bibitem{WeiDurh78} {\sc Wei, L. J.} and {\sc Durham,} S. (1978). The
randomized pay-the-winner rule in medical trials. {\em J. Amer.
Statist. Assoc.} {\bf 73} 840-843.

\bibitem{YW} {\sc Yao, Q. \& Wei, L. J.} (1996). Play the winner for
phase II/III clinical trials. {\em Statistics in Medicine} {\bf 15}:
2413-2423.

\bibitem{Zelen} {\sc Zelen, M.}  (1969).  Play-the-winner rule and
the controlled clinical trial.  {\em J. Amer. Statist. Assoc.} {\bf
64} 131-146.

\bibitem{Zhang04} {\sc Zhang, L.-X.} (2004). Strong approximations of martingale vectors and their
applications in Markov-chain adaptive designs. {\em Acta Math. Appl.
Sinica, English Series} {\bf 20}(2) 337--352.

\bibitem{Zhang06} {\sc  Zhang, L.-X.} (2006). Asymptotic results on a class of
adaptive multi-treatment designs. {\em J. Multi. Anal.} {\bf 97}
586--605.

\bibitem{ZCCH}  {\sc Zhang, L.-X., Chan, W.S., Cheung. S.H. \& Hu, F.} (2007). A
A generalized drop-the-loser urn for clinical trials with delayed
responses {\em Statistica Sinica} {\bf 15} in press.


\bibitem{ZHC} {\sc Zhang, L.-X., Hu, F., \& Cheung. S.H.} (2006). Asymptotic theorems of
sequential estimation-adjusted urn models. {\em Ann. Appl. Probab.}
{\bf 16}, 340-369.

\bibitem{ZHCC06a} {\sc Zhang, L-.X., Hu, F., Cheung, S.H. \& Chan,
W.S.} (2006a). Doubly adaptive biased coin designs with delayed
responses. {\em Manuscript}.

\bibitem{ZHCC06b} {\sc Zhang, L-.X., Hu, F., Cheung, S.H. \& Chan,
W.S.} (2006b). Immigrated urn models-- asymptotic properties and
applications. {\em Manuscript}.

\bibitem{ZHCC07} {\sc Zhang, L-.X., Hu, F., Cheung, S.H. \& Chan, W.S.}
 (2007). Asymptotic properties of covariate-adjusted response-adaptive designs. {\em Ann.
 Statist.} to appear.
 \underline{http://arxiv.org/abs/math.ST/0610518}


\end{thebibliography}
\end{document}